\documentclass[11pt]{amsart}
\usepackage{natbib}
\usepackage{amssymb}
\theoremstyle{definition}
\newtheorem{definition}{Definition}[section]
\newtheorem{example}[definition]{Example}
\newtheorem{remark}[definition]{Remark}
\theoremstyle{plain}
\newtheorem{theorem}[definition]{Theorem}
\newtheorem{lemma}[definition]{Lemma}
\newtheorem{corollary}[definition]{Corollary}
\newtheorem{proposition}[definition]{Proposition}
\newtheorem{conjecture}[definition]{Conjecture}

\newenvironment{veryspecialcolumnspacing}{\setlength\arraycolsep{1.2pt}}

\def\dispskip{\;\;\;}
\def\cc{(\mathcal C)}
\def\dc{(\mathcal D)}
\def\commutator#1#2{[#1,\,#2]}
\def\associator#1#2#3{[#1,\,#2,\,#3]}
\def\cdeg{\operatorname{cdeg}}
\def\rad{\operatorname{Rad}}
\def\to{\longrightarrow}
\def\aut{\operatorname{Aut}}
\def\chn#1{\overline{#1}}

\title[Moufang Loops that Share Associator]{Moufang Loops that Share Associator
and Three Quarters of their Multiplication Tables}
\author[A.~Dr\'apal]{Ale\v{s} Dr\'apal}
\email{drapal@karlin.mff.cuni.cz}
\address{Department of Algebra, Charles University, Sokolovsk\'a 83, 186 75 Prague,
Czech Republic}

\author[P.~Vojt\v{e}chovsk\'y]{Petr Vojt\v{e}chovsk\'y}
\email{petr@math.du.edu}
\address{Department of Mathematics, University of Denver, 2360 S Gaylord St,
Denver, CO 80208, USA}

\thanks{Work supported by Grant Agency of Charles University, grant number
$269/2001/$B-MAT/MFF. The first author supported also by institutional grant
MSM $113200007$.}

\begin{document}

\begin{abstract}
Two constructions due to Dr\'apal produce a group by modifying exactly one
quarter of the Cayley table of another group. We present these constructions in
a compact way, and generalize them to Moufang loops, using loop extensions.
Both constructions preserve associators, the associator subloop, and the
nucleus. We conjecture that two Moufang 2-loops of finite order $n$ with
equivalent associator can be connected by a series of constructions similar to
ours, and offer empirical evidence that this is so for $n=16$, $24$, $32$; the
only interesting cases with $n\le 32$. We further investigate the way the
constructions affect code loops and loops of type $M(G,\,2)$. The paper closes
with several conjectures and research questions concerning the distance of
Moufang loops, classification of small Moufang loops, and generalizations of
the two constructions. \vskip 2mm

\end{abstract}

\maketitle

\section{Introduction}

\noindent Moufang loops, i.e., loops satisfying the \emph{Moufang identity}
$((xy)x)z=x(y(xz))$, are surely the most extensively studied loops. Despite
this fact, the classification of Moufang loops is finished only for orders less
than $64$, and several ingenious constructions are needed to obtain all these
loops. The purpose of this paper is to initiate a new approach to finite
Moufang $2$-loops. Namely, we intend to decide whether all Moufang $2$-loops of
given order with equivalent associator can be obtained from just one of them,
using only group-theoretical constructions. (See below for details). We prove
that this is the case for $n=16$, $24$, and $32$, which are the only orders
$n\le 32$ for which there are at least two non-isomorphic nonassociative
Moufang loops ($5$, $5$, and $71$, respectively). We also show that for every
$m\ge 6$ there exist classes of loops of order $2^m$ that satisfy our
hypothesis. Each of these classes consists of code loops whose nucleus has
exactly two elements (cf.\ Theorem \ref{Th:Code}).

As it turns out, we will only need two constructions that were introduced in
\cite{Drapal}, and that we call \emph{cyclic} and \emph{dihedral}. They are
recalled in Sections \ref{Sc:Cyclic} and \ref{Sc:Dihedral}, and generalized to
Moufang loops in Sections \ref{Sc:MoufangCyclic} and \ref{Sc:MoufangDihedral}.
The main feature of both constructions is that, given a Moufang loop
$(G,\,\cdot)$, they produce a generally non-isomorphic Moufang loop $(G,\,*)$
that has the same associator and nucleus as $(G,\,\cdot)$, and whose
multiplication table agrees with the multiplication table of $(G,\,\cdot)$ in
$3/4$ of positions.

The constructions allow a very compact description with the help of simple
modular arithmetic, developed in Section \ref{Sc:Modular}. Nevertheless, in
order to prove that the constructions are meaningful for Moufang loops
(Theorems \ref{Th:Cyclic}, \ref{Th:Dihedral}), one benefits from knowing some
loop extension theory (Section \ref{Sc:FactorSets}). (An alternative proof
using only identities is available as well \cite{PV}, but is much longer.)

We then turn our attention to two classes of Moufang loops: code loops (Section
\ref{Sc:Code}), and loops of type $M(G,\,2)$ (Section \ref{Sc:Chein}).

Up to isomorphism, code loops can be identified with maps $P:V\to F$ whose
$3$rd derived form is trilinear, where $F=GF(2)$ and $V$ is a finite vector
space over $F$. Section \ref{Sc:Code} explains how $P$ is modified under our
constructions. These modification can be described in terms of linear and
quadratic forms, and it is not difficult to see how one can gradually transform
a code loop to any other code loop with equivalent associator (cf.\ Proposition
\ref{Pr:AllCode}).

The loops of type $M(G,\,2)$ play a prominent role in the classification of
Moufang loops, chiefly thanks to their abundance among small loops. In Section
\ref{Sc:Chein}, we describe how the loops $M(G,\,2)$ behave under both
constructions.

It has been conjectured \cite{D2} that from each finite $2$-group one can
obtain all other $2$-groups of the same order by repeatedly applying a
construction that preserves exactly $3/4$ of the corresponding multiplication
tables. For $n\le 32$, this conjecture is known to be true, and for such $n$ it
suffices to use only the cyclic and dihedral constructions \cite{Zhukavets}.
For $n=64$, these constructions yield two blocks of groups and it is not known
at this moment if there exists a similar construction that would connect these
two blocks \cite{BalekDrapal}.

In view of these results about $2$-groups, it was natural to ask how universal
the cyclic and dihedral constructions remain for Moufang loops of small order.
A computer search (cf.\ Section \ref{Sc:Small}) has shown that for orders
$n=16$, $24$, $32$ the blocks induced by cyclic and dihedral constructions
coincide with blocks of Moufang loops with equivalent associator. This is the
best possible result since none of the constructions changes the associator,
and since the two constructions are not sufficient even for groups when $n=64$.

The search for pairs of $2$-groups that can be placed at quarter distance (a
phrase expressing that $3/4$ of the multiplication tables coincide) stems from
the discovery that two $2$-groups which differ in less than a quarter of their
multiplication tables are isomorphic \cite{D2}. We conjecture that this
property remains true for Moufang $2$-loops. Additional conjectures, together
with suggestions for future work, can be found at the end of the paper.

We assume basic familiarity with calculations in nonassociative loops and in
Moufang loops in particular. The inexperienced reader should consult
\cite{Pflugfelder}.

A word about the notation. The dihedral group $\langle a,\,b;\;
a^n=b^2=1,\,aba=b\rangle$ of order $2n$ will be denoted by $D_{2n}$, although
some of the authors we cite use $D_n$; for instance \cite{GMR}. We count the
Klein $4$-group among dihedral groups, and denote it also by $V_4$. The
generalized quaternion group $\langle a,\,b;\;
a^{2^{n-1}}=1,\,a^{2^{n-2}}=b^2,\,bab^{-1}=a^{-1}\rangle$ of order $2^n$ will
be denoted by $Q_{2^n}$. We often write $ab$ instead of $a\cdot b$. In fact,
following the custom, we use ``$\cdot$'' to indicate the order in which
elements are multiplied. For example, $a\cdot bc$ stands for $a(bc) = a\cdot
(b\cdot c)$.

\section{Modular Arithmetic and the Function $\sigma$}\label{Sc:Modular}

\noindent Let $m$ be a positive integer and $M$ the set $\{-m+1$, $-m+2$,
$\dots$, $m-1$, $m\}$. Denote by $\oplus$ and $\ominus$ the addition and
subtraction modulo $2m$ in $M$, respectively. More precisely, define
$\sigma:\mathbb{Z}\to \{-1$, $0$, $1\}$ by
\begin{equation*}
    \sigma(i) = \left\{\begin{array}{ll}
        1,&i>m,\\
        0,&i\in M,\\
        -1,&i<1-m,
    \end{array}\right.
\end{equation*}
and let
\begin{displaymath}
    i\oplus j = i+j-2m\sigma(i+j),\dispskip
    i\ominus j= i-j-2m\sigma(i-j),
\end{displaymath}
for any $i$, $j\in M$. In order to eliminate parentheses, we postulate that
$\oplus$ and $\ominus$ are more binding than $+$ and $-$. Observe that $1-i$
belongs to $M$ whenever $i$ does, and that $\sigma(1-i) = -\sigma(i)$.

We will need the following identities for $\sigma$ in Sections
$\ref{Sc:Cyclic}$ and $\ref{Sc:Dihedral}$:
\begin{eqnarray}
    \sigma(i+j)+\sigma(i\oplus j+k) &=&
        \sigma(j+k)+\sigma(i+j\oplus k),
        \label{Eq:AssocCyclicSigma}\\
  -\sigma(i+j)+\sigma(1-i\oplus j+k) &=&
       \sigma(1-j+k)-\sigma(i+j\ominus k) .
        \label{Eq:AssocDihedralSigma}
\end{eqnarray}
The identity $(\ref{Eq:AssocCyclicSigma})$ follows immediately from $(i\oplus
j)\oplus k = i\oplus(j\oplus k)$. To establish $(\ref{Eq:AssocDihedralSigma})$,
consider $(i\oplus j)\ominus k=i\oplus(j\ominus k)$. This yields $-\sigma(i+j)
- \sigma(i\oplus j-k) = - \sigma(j-k) -\sigma(i+j\ominus k)$. Since
$-\sigma(i\oplus j-k) =\sigma(1-i\oplus j+k)$ and $-\sigma(j-k) =
\sigma(1-j+k)$, we are done.

\section{The Cyclic Construction}\label{Sc:Cyclic}

\noindent Let us start with the less technical of the two constructions---the
cyclic one. We will work in the more general setting of Moufang loops, and take
full advantage of the function $\sigma$ defined in Section \ref{Sc:Modular}.

Let $G$ be a Moufang loop. Recall that $Z(G)$, the \emph{center} of $G$,
consists of all elements that commute and associate with all elements of $G$.
In more detail, given $x$, $y$, $z\in G$, the \emph{commutator}
$\commutator{x}{y}$ of $x$, $y$ (resp.\ the \emph{associator}
$\associator{x}{y}{z}$ of $x$, $y$, $z$) is the unique element $w\in G$
satisfying $xy=yx\cdot w$ (resp.\ $(xy)z=x(yz)\cdot w$). When three elements of
a Moufang loop associate in some order, they associate in any order. Hence
$Z(G)=\{x\in G;\;\commutator{x}{y}=\associator{x}{y}{z}=1$ for every $y$, $z\in
G\}$.

We say that $(G,\,S,\alpha,\,h)$ satisfies condition $\cc$ if
\begin{enumerate}
\item[-] $G$ is a Moufang loop,

\item[-] $S\unlhd G$, and $G/S=\langle \alpha \rangle$ is a cyclic group of
order $2m$,

\item[-] $1\ne h\in S\cap Z(G)$.
\end{enumerate}
Then we can view $G$ as the disjoint union $\bigcup_{i\in M}\alpha^i$, and
define a new multiplication $*$ on $G$ by
\begin{equation}\label{Eq:CyclicConstruction}
    x*y = xyh^{\sigma(i+j)},
\end{equation}
where $x\in\alpha^i$, $y\in\alpha^j$, and $i$, $j\in M$.

The resulting loop (that is Moufang, as we shall see) will be denoted by
$(G,\,*)$. Whenever we say that $(G,\,S,\,\alpha,\,h)$ satisfies $\cc$, we
assume that $(G,\,*)$ is defined by $(\ref{Eq:CyclicConstruction})$.

The following Proposition is a special case of Theorem \ref{Th:Cyclic}. We
present it here because the associative case is much simpler than the Moufang
case.

\begin{proposition}\label{Pr:CyclicConstruction}
When $G$ is a group and $(G,\,S,\,\alpha,\,h)$ satisfies $\cc$ then $(G,\,*)$
is a group.
\end{proposition}
\begin{proof}
Let $x\in\alpha^i$, $y\in\alpha^j$, $z\in\alpha^k$, for some $i$, $j$, $k\in
M$. Since $h\in Z(G)$, we have
\begin{equation}\label{Eq:CyclicX}
\begin{array}{rcl}
    (x*y)*z &=& (xy)z\cdot h^{\sigma(i+j)+\sigma(i\oplus j+k)},\\
    x*(y*z) &=& x(yz)\cdot h^{\sigma(j+k)+\sigma(i+j\oplus k)}.
\end{array}
\end{equation}
This follows from $(\ref{Eq:CyclicConstruction})$ and from the fact that
$xy\in\alpha^{i\oplus j}$, $yz\in\alpha^{j\oplus k}$. By
$(\ref{Eq:AssocCyclicSigma})$, $(G,\,*)$ is associative.
\end{proof}

\section{The Dihedral Construction}\label{Sc:Dihedral}

\noindent We proceed to the dihedral construction. Let $G$ be a Moufang loop,
and let $N(G)$ be the nucleus of $G$. Recall that $N(G)=\{x\in G;\;
\associator{x}{y}{z}=1$ for every $y$, $z\in G\}$, and that
$\associator{x}{y}{z}=1$ implies $\associator{y}{x}{z}=\associator{x}{z}{y}=1$
for every $x$, $y$, $z\in G$.

We say that $(G,\,S,\,\beta,\,\gamma,\,h)$ satisfies condition $\dc$ if
\begin{enumerate}
\item[-] $G$ is a Moufang loop,

\item[-] $S\unlhd G$ and $G/S$ is a dihedral group of order $4m$ (where we
allow $m=1$),

\item[-] $\beta$, $\gamma$ are involutions of $G/S$ such that
$\alpha=\beta\gamma$ is of order $2m$,

\item[-] $1\ne h\in S\cap Z(G_0)\cap N(G)$ and $hxh=x$ for some (and hence
every) $x\in G_1$, where $G_0=\bigcup_{i\in M}\alpha^i$, $G_1=G\setminus G_0$.
\end{enumerate}
We can then choose $e\in\beta$ and $f\in\gamma$, view $G$ as the disjoint union
$\bigcup_{i\in M}(\alpha^i\cup e\alpha^i)$ or $\bigcup_{j\in M}(\alpha^j\cup
\alpha^jf)$, and define a new multiplication $*$ on $G$ by
\begin{equation}\label{Eq:DihedralConstruction}
    x*y=xyh^{(-1)^r\sigma(i+j)},
\end{equation}
where $x\in\alpha^i\cup e\alpha^i$, $y\in(\alpha^j\cup\alpha^jf)\cap G_r$, $i$,
$j\in M$, and $r\in\{0,\,1\}$.

The resulting loop (again always Moufang) will be denoted by $(G,\,*)$. As in
the cyclic case, whenever we say that $(G,\,S,\,\beta,\,\gamma,\,h)$ satisfies
$\dc$, we assume that $(G,\,*)$ is defined by
$(\ref{Eq:DihedralConstruction})$.

Note that $*$ does not depend on the choice of $e\in\beta$ and $f\in\gamma$.
Also note that when $(G,\,S,\,\beta,\,\gamma,\,h)$ satisfies $\dc$, then
$(G_0,\,S,\,\alpha=\beta\gamma,\,h)$ satisfies $\cc$.

Since $G/S$ is dihedral, $\alpha$, $\beta$ and $\gamma$ satisfy
\begin{equation*}
    \beta\alpha^i=\alpha^{\ominus i}\beta,\,
    \gamma\alpha^i=\alpha^{\ominus i}\gamma,\,
    \beta\alpha^i=\alpha^{1-i}\gamma,\,
    \alpha^i\gamma = \beta\alpha^{1-i},
\end{equation*}
for any $i\in M$, where we write $\ominus i$ rather than $-i$ to make sure that
the exponents remain in $M$.

\begin{remark}\label{Rm:Both}
Although $\alpha$, $G_0$, $G_1$, $e$ and $f$ are not explicitly mentioned in
condition $\dc$, we will often refer to them. Strictly speaking, we did not
need to include $S$ among the parameters of any of the constructions, as it can
always be calculated from the remaining parameters. Finally, we will sometimes
find ourselves in a situation when we do not want to treat $\cc$ and $\dc$
separately. Let us therefore agree that $G_0=G_1=G$, $e=f=1$, and that $\beta$,
$\gamma$ are meaningless when $\cc$ applies.
\end{remark}

\begin{lemma}\label{Lm:EF}
Assume that $(G,\,S,\,\beta,\,\gamma,\,h)$ satisfies $\dc$. Then
$(ex)*y=e(x*y)$ and $(x*y)f=x*(yf)$ whenever $y\in N(G)$.
\end{lemma}
\begin{proof}
Choose $x\in\alpha^i\cup e\alpha^i$, $y\in(\alpha^j\cup \alpha^jf)\cap G_r$,
and note that $ex$ belongs to $\alpha^i\cup e\alpha^i$, while $yf$ belongs to
$(\alpha^j\cup \alpha^jf)\cap G_{r+1}$. For the sake of brevity, set
$t=h^{(-1)^r\sigma(i+j)}$. Then $(ex)*y = (ex)y\cdot t=e(xy)\cdot t = e(xy\cdot
t) = e(x*y)$, and $(x*y)f = (xy\cdot t)f = xy\cdot tf = xy\cdot ft^{-1} =
(xy)f\cdot t^{-1} = x(yf)\cdot t^{-1} = x*(yf)$, where we used $y\in N(G)$ and
$h\in N(G)$ several times.
\end{proof}

Similarly as in the cyclic case, Proposition \ref{Pr:DihedralConstruction} is a
special case of Theorem \ref{Th:Dihedral}:

\begin{proposition}\label{Pr:DihedralConstruction}
When  $G$ is a group and $(G,\,S,\,\beta,\,\gamma,\,h)$ satisfies $\dc$ then
$(G,\,*)$ is a group.
\end{proposition}

\begin{proof}
If $(x*y)*z=x*(y*z)$, Lemma \ref{Lm:EF} implies that $((ex)*y)*z=(ex)*(y*z)$
and $(x*y)*(zf)=x*(y*(zf))$. We can therefore assume that $x\in\alpha^i$,
$z\in\alpha^k$, and $y\in\alpha^j\cup \alpha^jf$, for some $i$, $j$, $k\in M$.

When $y\in\alpha^j$, the definition $(\ref{Eq:DihedralConstruction})$ of $*$
coincides with the cyclic case $(\ref{Eq:CyclicConstruction})$, and $x$, $y$,
$z$ associate in $(G,\,*)$ by Proposition \ref{Pr:CyclicConstruction}. Assume
that $y\in\alpha^jf\subseteq G_1$, and recall the coset relations
$\alpha^j\gamma=\beta\alpha^{1-j}$. Then
\begin{equation}\label{Eq:DihedralX}
\begin{array}{rcl}
    (x*y)*z &=& (xy)z\cdot h^{-\sigma(i+j)+\sigma(1-i\oplus j+k)},\\
    x*(y*z) &=& x(yz)\cdot h^{\sigma(1-j+k)-\sigma(i+j\ominus k)},
\end{array}
\end{equation}
because $xy\in\alpha^i\alpha^j\gamma=\alpha^{i\oplus
j}\gamma=\beta\alpha^{1-i\oplus j}$, and
$yz\in\alpha^j\gamma\alpha^k=\alpha^{j\ominus k}\gamma$. By
$(\ref{Eq:AssocDihedralSigma})$, $(G,\,*)$ is associative.
\end{proof}

\section{Factor Sets}\label{Sc:FactorSets}

\noindent Before we prove that $(G,\,*)$ is a Moufang loop if $\cc$ or $\dc$ is
satisfied, let us briefly review extensions of abelian groups by Moufang loops.
We follow closely the group-theoretical approach, cf.\ \cite[Ch.\
11]{Robinson}.

Let $Q$ be a Moufang loop and $A$ a $Q$-module. Since, later on, we will deal
with two extensions at the same time, we shall give a name to the action of $Q$
on $A$, say $\varphi: Q\to \aut A$. Consider a map $\eta:Q\times Q\to A$, and
define a new multiplication on the set product $Q\times A$ by
\begin{displaymath}\label{Eq:Extension}
    (x,\,a)(y,\,b) = (xy,\,a^{\varphi(y)}+ b +\eta(x,\,y)),
\end{displaymath}
where we use additive notation for the abelian group $A$. The resulting
quasigroup will be denoted by $E=(Q,\,A,\,\varphi,\,\eta)$.

It is easy to see that $E$ is a loop if and only if there exists $c\in A$ such
that
\begin{equation}\label{Eq:IsLoop}
    \eta(x,\,1) = c, \dispskip \eta(1,\,x) = c^{\varphi(x)},
\end{equation}
for every $x\in Q$. The neutral element of $E$ is then $(1,\,-c)$.

From now on, we will assume that $E$ satisfies $(\ref{Eq:IsLoop})$ with $c=0$,
and speak of $E$ as an \emph{extension} of $A$ by $Q$. Verify that $E$ is a
group if and only if $Q$ is a group and
\begin{equation}\label{Eq:IsGroup}
    \eta(x,\,y)^{\varphi(z)}+\eta(xy,\,z)=\eta(y,\,z)+\eta(x,\,yz)
\end{equation}
holds for every $x$, $y$, $z\in Q$. Moreover, using the Moufang identity
$(xy\cdot x)z = x(y\cdot xz)$, one can check by straightforward calculation
that $E$ is a Moufang loop if and only if
\begin{equation}\label{Eq:IsMoufang}
    \eta(x,\,y)^{\varphi(xz)}+\eta(xy,\,x)^{\varphi(z)}+\eta(xy\cdot x,\,z)
    =\eta(x,\,z)+\eta(y,\,xz)+\eta(x,\,y\cdot xz)
\end{equation}
holds for every $x$, $y$, $z\in Q$. (Note that $\varphi(y\cdot
xz)=\varphi(yx\cdot z)$ even if $x$, $y$, $z$ do not associate.)

Every pair $(\varphi,\,\eta)$ satisfying $(\ref{Eq:IsLoop})$ with $c=0$ is
called a \emph{factor set}. If it also satisfies $(\ref{Eq:IsGroup})$, resp.\
$(\ref{Eq:IsMoufang})$, we call it \emph{associative factor set}, resp.\
\emph{Moufang factor set}.

Given two factor sets $(\varphi$, $\eta)$ and $(\varphi$, $\mu)$, we can obtain
another factor set, their \emph{sum} $(\varphi$, $\eta+\mu)$, by letting
$(\eta+\mu)(x,\,y) = \eta(x,\,y)+\mu(x,\,y)$ for every $x$, $y\in Q$. Since $A$
is an abelian group, the sum of two associative factor sets (resp.\ Moufang
factor sets) is associative (resp.\ Moufang). As every group is a Moufang loop,
it must be the case that every associative factor set is Moufang. Here is a
proof that only refers to factor sets:

\begin{lemma}
Every associative factor set is Moufang.
\end{lemma}
\begin{proof}
Let $(\varphi$, $\eta)$ be an associative factor set. Substituting $xz$ for $z$
in $(\ref{Eq:IsGroup})$ yields
\begin{equation}\label{Eq:IsGroup1}
    \eta(x,\,y)^{\varphi(xz)}+\eta(xy,\,xz)=\eta(y,\,xz)+\eta(x,\,y\cdot xz),
\end{equation}
while substituting $xy$ for $x$, and simultaneously $x$ for $y$ in
$(\ref{Eq:IsGroup})$ yields
\begin{equation}\label{Eq:IsGroup2}
    \eta(xy,\,x)^{\varphi(z)}+\eta(xy\cdot x,\,z)=\eta(x,\,z)+\eta(xy,\,xz).
\end{equation}
The identity $(\ref{Eq:IsMoufang})$ is obtained by adding $(\ref{Eq:IsGroup1})$
to $(\ref{Eq:IsGroup2})$ and subtracting $\eta(xy,\,xz)$ from both sides.
\end{proof}

Assume that $(\varphi,\,\eta)$ is a Moufang factor set. Then the right inverse
of $(x,\,a)$ in $(Q,\,A,\,\varphi,\,\eta)$ is
$(x^{-1},\,-a^{\varphi(x^{-1})}-\eta(x,\,x^{-1}))$, as a short calculation
reveals. Similarly, the left inverse of $(x,\,a)$ is
$(x^{-1},\,-a^{\varphi(x^{-1})}-\eta(x^{-1},\,x)^{\varphi(x^{-1})})$. Since
$(Q,\,A,\,\varphi,\,\eta)$ is a Moufang loop, the two inverses coincide, and we
have
\begin{equation}\label{Eq:Inverses}
    \eta(x,\,x^{-1})=\eta(x^{-1},\,x)^{\varphi(x^{-1})},
\end{equation}
for any Moufang factor set $(\varphi,\,\eta)$ and $x\in Q$.
(Alternatively---and more naturally---the identity $(\ref{Eq:Inverses})$
follows immediately from $(\ref{Eq:IsMoufang})$ when we substitute $x^{-1}$ for
$x$, $x$ for $y$, and $1$ for $z$.)

\begin{lemma}\label{Lm:FactorAssociators}
Assume that $(\varphi,\,\eta)$ is a Moufang factor set and $(\varphi,\,\mu)$ is
an associative factor set. Then the associators in $(Q,\,A,\,\varphi,\,\eta)$
and $(Q,\,A,\,\varphi,\,\eta+\mu)$ coincide if and only if
\begin{equation}\label{Eq:ConditionOnMu}
    \mu((x\cdot yz)^{-1},\,xy\cdot z) = \mu(x\cdot yz,\,(x\cdot
    yz)^{-1})^{\varphi(xy\cdot z)}
\end{equation}
for every $x$, $y$, $z\in Q$. This happens if and only if
\begin{equation}\label{Eq:ShorterConditionOnMu}
    \mu(x\cdot yz,\,\associator{x}{y}{z})=0
\end{equation}
for every $x$, $y$, $z\in Q$. In particular, the associators coincide if $Q$ is
a group.
\end{lemma}

\begin{proof}
Let $(x,\,a)$, $(y,\,b)$, $(z,\,c)\in (Q,\,A,\,\varphi,\,\eta)$. Then
\begin{eqnarray*}
    u=(x,\,a)(y,\,b)\cdot (z,\,c) &=& (xy\cdot z,\,s+t), \\
    v=(x,\,a)\cdot (y,\,b)(z,\,c) &=& (x\cdot yz,\,s),
\end{eqnarray*}
where
\begin{eqnarray*}
    s&=&a^{\varphi(yz)}+b^{\varphi(z)}+c+\eta(y,\,z)+\eta(x,\,yz),\\
    t&=&\eta(x,\,y)^{\varphi(z)} + \eta(xy,\,z)-\eta(y,\,z)-\eta(x,\,yz).
\end{eqnarray*}

The associator $\associator{(x,\,a)}{(y,\,b)}{(z,\,c)}$ in
$(Q,\,A,\,\varphi,\,\eta)$ is therefore equal to $v^{-1}u =
(\associator{x}{y}{z},\,d)$, where
\begin{displaymath}
    d = t+\eta((x\cdot yz)^{-1},\,xy\cdot z) -\eta(x\cdot yz,\,(x\cdot
    yz)^{-1})^{\varphi(xy\cdot z)}.
\end{displaymath}
 Similarly, the same associator in
$(Q,\,A,\,\varphi,\,\eta+\mu)$ is $(\associator{x}{y}{z},\,d+e+f)$, where
\begin{eqnarray*}
    e&=&\mu(x,\,y)^{\varphi(z)}+\mu(xy,\,z)-\mu(y,\,z)-\mu(x,\,yz),\\
    f&=&\mu((x\cdot yz)^{-1},\,xy\cdot z)-\mu(x\cdot yz,\,(x\cdot
    yz)^{-1})^{\varphi(xy\cdot z)}.
\end{eqnarray*}
Since $(\varphi,\,\mu)$ satisfies $(\ref{Eq:IsGroup})$, $e$ vanishes. Therefore
the two associators coincide for all $x$, $y$, $z\in Q$ if and only if
$(\ref{Eq:ConditionOnMu})$ is satisfied for every $x$, $y$, $z\in Q$.

Substituting $x\cdot yz$ for $x$, $(x\cdot yz)^{-1}$ for $y$, and $xy\cdot z$
for $z$ into $(\ref{Eq:IsGroup})$ yields
\begin{displaymath}
    \mu(x\cdot yz,\,(x\cdot yz)^{-1})^{\varphi(xy\cdot z)}
    = \mu((x\cdot yz)^{-1},\,xy\cdot z)
    + \mu(x\cdot yz,\,\associator{x}{y}{z}).
\end{displaymath}
Hence $(\ref{Eq:ConditionOnMu})$ is satisfied if and only if
$(\ref{Eq:ShorterConditionOnMu})$ holds. The latter condition is of course
satisfied when $Q$ is a group.
\end{proof}

\section{The Cyclic Construction for Moufang loops}\label{Sc:MoufangCyclic}

\noindent Throughout this section, assume that $(G,\,S,\,\alpha,\,h)$ satisfies
$\cc$, and that $A$ is the subloop of $S$ generated by $h$. Using loop
extensions, we prove that $(G,\,*)$ is a Moufang loop with the same
associators, associator subloop, and nucleus as $(G,\,\cdot)$. Recall that the
\emph{associator subloop} of a loop $L$ is the subloop $A(L)$ generated by all
associators $\associator{x}{y}{z}$, where $x$, $y$, $z\in L$.

\begin{lemma}\label{Lm:CycNormal}
$A$ is a normal subloop of both $(G,\,\cdot)$ and $(G,\,*)$. Moreover,
$(G,\,\cdot)/A=(G,\,*)/A$.
\end{lemma}
\begin{proof}
Since $h\in Z(G,\,\cdot)$, the subgroup $A=\langle h \rangle \subseteq
Z(G,\,\cdot)$ is normal in $(G,\,\cdot)$. In fact, $x*h=xh$, $h*x=hx$ for every
$x\in G$ (since $h\in S=\alpha^0$), and thus $A$ is normal in $(G,\,*)$ as
well.

Write the elements of $G/A$ as cosets $xA$. Since, for some $t$, we have
$xA\cdot yA = (xy)A$ and $xA*yA = (x*y)A=(xyh^t)A = (xy)A$, the loops
$(G,\,\cdot)/A$ and $(G,\,*)/A$ coincide.
\end{proof}

Let $Q$ be the Moufang loop $(G,\,\cdot)/A=(G,\,*)/A$. Let $\iota$ be the
trivial homomorphism $Q\to \aut A$, $\iota(q)=id_A$, for every $q\in Q$. We
want to construct two factor sets $(\iota,\,\eta)$, $(\iota,\,\eta^*)$ such
that $(Q,\,A,\,\iota,\,\eta)\simeq (G,\,\cdot)$ and
$(Q,\,A,\,\iota,\,\eta^*)\simeq (G,\,*)$. In order to save space, we keep
writing the operation in $A$ multiplicatively.

Let $\pi:Q=G/A\to G$ be a transversal, i.e., a map satisfying $\pi(xA)\in xA$
for every $x\in G$. Then, for every $xA$, $yA$, there is an integer
$\tau(xA,\,yA)$ such that $\pi((xy)A) = \pi(xA)\pi(yA)h^{\tau(xA,\,yA)}$.

\begin{proposition}\label{Pr:CyclicFactor}
Assume that $(G,\,S,\,\alpha,\,h)$ satisfies $\cc$, and that $A$ is the subloop
of $S$ generated by $h$. With $Q=(G,\,\cdot)/A=(G,\,*)/A$ and $\tau$ as above,
define $\eta$, $\eta^*:Q\times Q\to A$ by
\begin{eqnarray*}
    \eta(xA,\,yA) &=& h^{-\tau(xA,\,yA)},\\
    \eta^*(xA,\,yA) &=& \eta(xA,\,yA)h^{\sigma(i+j)},
\end{eqnarray*}
where $x\in\alpha^i$, $y\in\alpha^j$, and $i$, $j\in M$. Then
$(Q,\,A,\,\iota,\,\eta)\simeq (G,\,\cdot)$ and $(Q,\,A,\,\iota,\,\eta^*)\simeq
(G,\,*)$.
\end{proposition}
\begin{proof}
First of all, when $x$ belongs to $\alpha^i$ then every element of $xA$ belongs
to $\alpha^i$, and so $\eta^*$ is well-defined.

Let $\theta:(Q,\,A,\,\iota,\,\eta)\to (G,\,\cdot)$ be defined by
$\theta(xA,\,h^a) = \pi(xA)h^a$. Note that $\theta$ is well-defined, and that
it is clearly a bijection. Since
\begin{eqnarray*}
    &&\theta((xA,\,h^a)(yA,\,h^b))= \theta((xy)A,\,h^{a+b}\eta(xA,\,yA))
        =\pi((xy)A)h^{a+b}\eta(xA,\,yA)\\
    &&= \pi(xA)\pi(yA)h^{\tau(xA,\,yA)}h^{a+b}h^{-\tau(xA,\,yA)}
      = \pi(xA)h^a\pi(yA)h^b = \theta(xA,\,h^a)\theta(yA,\,h^b),
\end{eqnarray*}
$\theta$ is an isomorphism.

Similarly, let $\theta^*:(Q,\,A,\,\iota,\,\eta^*)\to (G,\,*)$ be defined by
$\theta^*(xA,\,h^a)=\pi(xA)h^a$. This is again a bijection. Pick
$x\in\alpha^i$, $y\in\alpha^j$. Since
\begin{eqnarray*}
    &&\theta^*((xA,\,h^a)(yA,\,h^b)) = \theta^*((xy)A,\,h^{a+b}\eta^*(xA,\,yA))
        = \pi((xy)A)h^{a+b}\eta^*(xA,\,yA)\\
    &&\phantom{mmm}=\pi(xA)\pi(yA)h^{\tau(xA,\,yA)}h^{a+b}
        h^{-\tau(xA,\,yA)}h^{\sigma(i+j)}
        = \pi(xA)\pi(yA)h^{a+b}h^{\sigma(i+j)}\\
    &&\phantom{mmm}= \pi(xA)h^a*\pi(yA)h^b = \theta^*(xA,\,h^a)*\theta^*(yA,\,h^b),
\end{eqnarray*}
$\theta^*$ is an isomorphism.
\end{proof}

We are now ready to prove the main theorem for the cyclic construction:

\begin{theorem}\label{Th:Cyclic}
The Moufang factor sets $(\iota$, $\eta)$ and $(\iota$, $\eta^*)$ introduced in
Proposition $\ref{Pr:CyclicFactor}$ differ by an associative factor set
$(\iota,\,\mu)$ that satisfies $(\ref{Eq:ConditionOnMu})$. Consequently,
$(G,\,*)$ is a Moufang loop, the associators in $(G,\,\cdot)$ and $(G,\,*)$
coincide, $A(G,\,\cdot) = A(G,\,*)$ coincide as loops, and $N(G,\,\cdot) =
N(G,\,*)$ coincide as sets.
\end{theorem}
\begin{proof}
With $\mu=\eta^*-\eta$ and $x\in\alpha^i$, $y\in\alpha^j$, we have
$\mu(xA,\,yA)=h^{\sigma(i+j)}$. Since $\mu(xA,\,A) = \mu(A,\,xA) =
h^{\sigma(i)} = h^0 = 1$, $(\iota,\,\mu)$ is a factor set. Pick further
$z\in\alpha^k$. We must verify that $(\iota,\,\mu)$ is associative, i.e., that
\begin{displaymath}
    \mu(xA,\,yA)\mu(xA,\,yAzA) = \mu(yA,\,zA)\mu(xA,\,yAzA).
\end{displaymath}
But this follows immediately from $(\ref{Eq:AssocCyclicSigma})$, as
$xAyA\in\alpha^{i\oplus j}$ and $yAzA\in\alpha^{j\oplus k}$. Thus
$(\iota,\,\mu)$ is associative, in particular Moufang. Then $(\iota$,
$\eta^*)=(\iota$, $\eta) + (\iota$, $\mu)$ is a Moufang factor set.

It is easy to verify that all associators of $(G,\,\cdot)$ belong to
$\alpha^0$. This means that $\mu(xAyA\cdot zA,\,\associator{xA}{yA}{zA})$
vanishes, and hence the associators in $(G,\,\cdot)$ and $(G,\,*)$ coincide by
Lemma \ref{Lm:FactorAssociators}. The associator subloops $A(G,\,\cdot)$ and
$A(G,\,*)$ are therefore generated by the same elements. In fact, the
multiplication in $A(G,\,\cdot)$ coincides with the multiplication in
$A(G,\,*)$ because, once again, every associator belongs to $\alpha^0$.
Finally, since an element belongs to the nucleus if and only if it associates
with all other elements, we must have $N(G,\,\cdot)=N(G,\,*)$.
\end{proof}

\section{The Dihedral Construction for Moufang Loops}\label{Sc:MoufangDihedral}

\noindent We are now going to prove that the dihedral construction works for
Moufang loops, too. The reasoning is essentially that of Section
\ref{Sc:MoufangCyclic}, however, we decided that it deserves a separate
treatment since it differs in several details. The confident reader can proceed
directly to the next section.

Throughout this section, we assume that $(G,\,S,\,\,\beta,\,\gamma,\,h)$
satisfies $\dc$, and that $A$ is the subloop of $S$ generated by $h$.

\begin{lemma}\label{Lm:DihNormal}
$A$ is a normal subloop of both $(G,\,\cdot)$ and $(G,\,*)$. Moreover,
$(G,\,\cdot)/A =(G,\,*)/A$.
\end{lemma}
\begin{proof}
We claim that $A$ is a normal subloop of $(G,\,\cdot)$. It suffices to prove
that $xA=Ax$, $x(Ay)=(xA)y$ and $x(yA)=(xy)A$ for every $x$, $y\in G$. Since
$A\le N(G)$, we only have to show that $xA=Ax$ for every $x\in G$. When $x\in
G_0$, there is nothing to prove as $h\in Z(G_0)$. When $x\in G_1$, we have
$xA=\{xh^a;\; 0\le a<2m\} = \{h^{-a}x;\; 0\le a<2m\} = Ax$, because $hxh=x$.
Thus $A$ is normal in $(G,\,\cdot)$. In fact, $x*h=xh$, $h*x=hx$ for every
$x\in G$ (since $h\in S=\alpha^0$), and thus $A$ is normal in $(G,\,*)$ as
well.

Write the elements of $G/A$ as cosets $xA$. Since, for some $t$, we have
$xA\cdot yA=(xy)A$ and $xA*yA = (x*y)A = (xyh^t)A = (xy)A$, the loops
$(G,\,\cdot)/A$ and $(G,\,*)/A$ coincide.
\end{proof}

We let $Q$ be the Moufang loop $(G,\,\cdot)/A=(G,\,*)/A$, and continue to
construct two factor sets $(\varphi,\,\eta)$, $(\varphi,\,\eta^*)$ such that
$(Q,\,A,\,\varphi,\,\eta)\simeq (G,\,\cdot)$ and
$(Q,\,A,\,\varphi,\,\eta^*)\simeq (G,\,*)$.

Fix a transversal $\pi:Q=G/A\to G$. Then, for every $xA$, $yA$, there is an
integer $\tau(xA,\,yA)$ such that $\pi((xy)A) =
\pi(xA)\pi(yA)h^{\tau(xA,\,yA)}$.

\begin{proposition}\label{Pr:DihedralFactor}
Assume that $(G,\,S,\,\beta,\,\gamma,\,h)$ satisfies $\dc$, and that $A$ is the
subloop of $S$ generated by $h$. With $Q=(G,\,\cdot)/A=(G,\,*)/A$ and $\tau$ as
above, define $\varphi:Q\to \aut A$ by $a^{\varphi(y)} = a^{(-1)^r}$, where
$y\in G_r$, $r\in\{0,\,1\}$. Furthermore, define $\eta$, $\eta^*:Q\times Q\to
A$ by
\begin{eqnarray*}
    \eta(xA,\,yA) &=& h^{-\tau(xA,\,yA)},\\
    \eta^*(xA,\,yA) &=& \eta(xA,\,yA)h^{(-1)^r\sigma(i+j)},
\end{eqnarray*}
where $x\in\alpha^i\cup e\alpha^i$, $y\in(\alpha^j\cup \alpha^jf)\cap G_r$,
$i$, $j\in M$, $r\in\{0,\,1\}$.  Then $(Q,\,A,\,\varphi,\,\eta)\simeq
(G,\,\cdot)$ and $(Q,\,A,\,\varphi,\,\eta^*)\simeq (G,\,*)$.
\end{proposition}
\begin{proof}
Since $G_rG_s=G_{r+s\pmod 2}$ for every $r$, $s\in\{0,\,1\}$, $\varphi$ is a
homomorphism.

When $x$ belongs to $\alpha^i\cup e\alpha^i$, then every element of $xA$
belongs to $\alpha^i\cup e\alpha^i$. When $y$ belongs to $(\alpha^j\cup
\alpha^jf)\cap G_r$, then every element of $yA$ belongs to
$(\alpha^j\cup\alpha^jf)\cap G_r$. Hence $\eta^*$ is well-defined.

Let $\theta:(Q,\,A,\,\varphi,\,\eta)\to (G,\,\cdot)$ be defined by
$\theta(xA,\,h^a)=\pi(xA)h^a$. This is clearly a well-defined bijection. When
$y\in G_r$, we have
\begin{eqnarray*}
    &&\theta((xA,\,h^a)(yA,\,h^b))
        = \theta((xy)A,\,h^{(-1)^r a}h^b\eta(xA,\,yA))\\
    &&= \pi((xy)A)h^{(-1)^r a}h^b\eta(xA,\,yA)
        = \pi(xA)\pi(yA)h^{\tau(xA,\,yA)}h^{(-1)^r a}h^bh^{-\tau(xA,\,yA)}\\
    &&= \pi(xA)\pi(yA)h^{(-1)^ra}h^b
        = \pi(xA)h^a\pi(yA)h^b = \theta(xA,\,h^a)\theta(yA,\,h^b),
\end{eqnarray*}
and $\theta$ is an isomorphism.

Similarly, let $\theta^*:(Q,\,A,\,\varphi,\,\eta^*)\to (G,\,*)$ be defined by
$\theta^*(xA,\,h^a)=\pi(xA)h^a$. This is again a bijection. With
$x\in\alpha^i\cup e\alpha^i$, $y\in(\alpha^j\cup\alpha^jf)\cap G_r$, we have
\begin{eqnarray*}
    &&\theta^*((xA,\,h^a)(yA,\,h^b))
        =\theta^*((xy)A,\,h^{(-1)^r a}h^b\eta^*(xA,\,yA))\\
    &&=\pi((xy)A)h^{(-1)^r a}h^b\eta^*(xA,\,yA)\\
    &&=\pi(xA)\pi(yA)h^{\tau(xA,\,yA)}h^{(-1)^ra}h^b
        h^{-\tau(xA,\,yA)}h^{(-1)^r \sigma(i+j)}\\
    &&=\pi(xA)h^a\pi(yA)h^bh^{(-1)^r \sigma(i+j)}
      = \pi(xA)h^a*\pi(yA)h^b = \theta^*(xA,\,h^a)*\theta^*(yA,\,h^b),
\end{eqnarray*}
and $\theta^*$ is an isomorphism.
\end{proof}

\begin{theorem}\label{Th:Dihedral}
The Moufang factor sets $(\varphi,\,\eta)$ and $(\varphi,\,\eta^*)$ introduced
in Proposition $\ref{Pr:DihedralFactor}$ differ by an associative factor set
$(\varphi,\,\mu)$ that satisfies $(\ref{Eq:ConditionOnMu})$. Consequently,
$(G,\,*)$ is a Moufang loop, the associators in $(G,\,\cdot)$ and $(G,\,*)$
coincide, $A(G,\,\cdot) = A(G,\,*)$ coincide as loops, and $N(G,\,\cdot) =
N(G,\,*)$ coincide as sets.
\end{theorem}
\begin{proof}
Let $\mu=\eta^*-\eta$. For $x\in\alpha^i\cup e\alpha^i$, $y\in(\alpha^j\cup
\alpha^jf)\cap G_r$, we have $\mu(xA,\,yA)=h^{(-1)^r \sigma(i+j)}$.

Since $\mu(xA,\,A)=\mu(A,\,xA)=h^0=1$, $(\varphi,\,\mu)$ is a factor set. By
the first 2 paragraphs of the proof of Proposition
\ref{Pr:DihedralConstruction}, $(\varphi,\,\mu)$ is associative, hence Moufang.
Then $(\varphi,\,\eta^*) = (\varphi,\,\eta)+(\varphi,\,\mu)$ is a Moufang
factor set.

It is easy to verify that every associator of $(G,\,\cdot)$ belongs to
$\alpha^0$. We can therefore reach the same conclusion as in Theorem
\ref{Th:Cyclic}.
\end{proof}

\section{Code Loops}\label{Sc:Code}

\noindent Now when we know that $(G,\,*)$ is a Moufang loop for both
constructions, we will focus on the effect the constructions have on two
important classes of Moufang loops: code loops and loops of type $M(G,\,2)$.
These loops are abundant among small Moufang loops, as we will see in Section
$\ref{Sc:Small}$. The results of Sections \ref{Sc:Code} and \ref{Sc:Chein} are
not needed elsewhere in this paper. Let us get started with code loops.

A loop $G$ is called \emph{symplectic} if it possesses a central subloop $Z$ of
order $2$ such that $G/Z$ is an elementary abelian $2$-group. When $G$ is
symplectic, we can define $P:G/Z\to Z$, $C:G/Z\times G/Z\to Z$, $A:G/Z\times
G/Z\times G/Z\to Z$ by $P(aZ)=a^2$, $C(aZ,\,bZ)=\commutator{a}{b}$,
$A(aZ,\,bZ,\,cZ)=\associator{a}{b}{c}$, for every $a$, $b$, $c\in G$. Note that
the three maps are well defined. For obvious reasons, we will often call $P$
the \emph{power map}, $C$ the \emph{commutator map}, and $A$ the
\emph{associator map}.

Every symplectic loop $G$ is an extension $(V,\,F,\,\iota,\,\eta)$ of the
$2$-element field $F=\{0,\,1\}$ by a finite vector space $V$ over $F$, where
$\eta:V\times V\to F$ satisfies $\eta(u,\,0)=\eta(0,\,u)=0$ for every $u\in V$
(i.e., $(\iota,\,\eta)$ is a factor set as defined in Section
\ref{Sc:FactorSets}). We can then identify $F$ with $Z$, $V$ with $G/Z$, and
consider $P$, $C$, $A$ as maps $P:V\to F$, $C:V\times V\to F$, $A:V\times
V\times V\to F$.

It is known that the triple $(P,\,C,\,A)$ determines the isomorphism type of
$G$ (cf.\ \cite[Theorem 12.13]{Aschbacher}).

Before we introduce code loops, we must define derived forms and combinatorial
degree. We will restrict the definitions to the two-element field $F$; more
general definitions can be found in \cite{Aschbacher} and \cite{Ward}.

Let $f:V\to F$ be a map satisfying $f(0)=0$. Then the \emph{$n$th derived form}
$f_n:V^n\to F$ of $f$ is defined by
\begin{displaymath}
    f_n(v_1,\,\dots,\,v_n)=
        \sum_{\{i_1,\,\dots,\,i_m\}\subseteq\{1,\,\dots,\,n\}}
        f(v_{i_1}+\cdots + v_{i_m}),
\end{displaymath}
where the summation runs over all nonempty subsets of $\{1,\,\dots,\,n\}$.
Although it is not immediately obvious, $f_n(v_1,\,\dots,\,v_n)$ vanishes
whenever $v_1$, $\dots$, $v_n$ are linearly dependant, and it makes sense to
define the \emph{combinatorial degree} of $f$, $\cdeg f$, as the smallest
nonnegative integer $n$ such that $f_{n+1}=0$.

Every form $f_n$ is symmetric, and two consecutive derived forms are related by
\emph{polarization}, i.e.,
\begin{displaymath}
    f_{n+1}(v_1,\,\dots,\,v_{n+1})=
        f_n(v_1,\,v_3,\,\dots,\,v_{n+1}) + f_n(v_2,\,\dots,\,v_{n+1})
        + f_n(v_1+v_2,\,v_3,\,\dots,\,v_{n+1}),
\end{displaymath}
for every $v_1$, $\dots$, $v_{n+1}\in V$. Thus $f_n$ is $n$-linear if and only
if $\cdeg{f}\le n$. Since $f(0)=0$, the form $f_2$ is alternating. Recall that
every alternating bilinear form over the two-element field is symmetric. When
$f$ is a quadratic form, $f_2$ is an alternating (thus symmetric) bilinear
form. Therefore the subspace of all forms $f:V\to F$ with $\cdeg{f}\le 2$
coincides with the subspace of all quadratic forms.

A symplectic loop $G$ defined on $V\times F$ is called a \emph{code loop} if
the power map $P:V\to F$ has $\cdeg P\le 3$, the commutator map $C$ coincides
with $P_2$, and the associator map $A$ coincides with $P_3$. The power map
therefore determines a code loop up to an isomorphism, and we will use the
notation $G=(V,F,P)$.

\begin{remark} Code loops were discovered by Griess \cite{Griess}, who used them to
elucidate the construction of the Parker loop, that is in turn involved in the
construction of the Monster group. We completely ignore the code aspect of code
loops here, and model our approach on \cite{Aschbacher} and \cite{Hsu}.
\end{remark}

Of course, not every symplectic loop is a code loop, however, as Aschbacher
proved in \cite[Lemma 13.1]{Aschbacher}, Chein and Goodaire in
\cite{CheinGoodaire}, and Hsu in \cite{Hsu}:

\begin{theorem}\label{Th:Hsu}
Code loops are exactly symplectic Moufang loops.
\end{theorem}

Thus our two constructions apply to code loops and we proceed to have a closer
look at them. Recall that the \emph{radical} $\rad{f}$ of an $n$-linear form
$f:V^n\to F$ is the subspace consisting of all vectors $v_1\in V$ such that
$f(v_1,\,\dots,\,v_n)=0$ for every $v_2$, $\dots$, $v_n\in V$.

The radical of $P_3$ determines the nucleus of the associated code loop, and
vice versa. We offer a complete description of the situation when $P_3$ has
trivial radical (i.e., $\rad{P_3}=F$). Then there is only one choice of $h$ for
$\cc$ and $\dc$ (see below). We expect to return to code loops with nontrivial
radical in a future paper.

\begin{remark}
Code loops with nontrivial radical are not closed under the two constructions.
$($cf.\ Example $\ref{Ex:Code})$. In fact, all code loops of order $32$ have
this property.
\end{remark}

\begin{lemma}\label{Lm:Code} Let $G=(V,F,P)$ be a code loop.
Assume that $\cc$ or $\dc$ is satisfied with some $h$, $S$. Then:
\begin{enumerate}
\item[(i)] If $G$ is not a group or if $h\in F$, then $S\supseteq F$, and
$G/S\simeq C_2$ or $G/S\simeq V_4$.

\item[(ii)] If $h\in F$ then the resulting loop $(G,\,*)$ is a code loop with
the same radical as $G$.

\item[(iii)] If $\rad{P_3}=F$ then $h\in N(G)=Z(G)=F$.
\end{enumerate}
\end{lemma}
\begin{proof}
Since $G=(V,F,P)$ is a code loop, we have $A(G)\subseteq F$. Let us prove (i).
First assume that $G$ is not a group. Since $|F|=2$, we must have $A(G)=F$. As
$G/S$ is associative, the subloop $S$ contains $A(G)=F$. Now assume that $1\ne
h\in F$. Since $h$ belongs to $S$, we immediately obtain $S\supseteq F$. Hence,
in any case, $G/S\le G/F$, and $G/S$ is an elementary abelian $2$-group. The
only two elementary abelian $2$-groups satisfying $\cc$ or $\dc$ are $C_2$ and
$V_4$, respectively.

To prove (ii), assume that $h\in F$. Then $(F,\,*)$ is a subloop of $(G,\,*)$,
by $(\ref{Eq:CyclicConstruction})$ and $(\ref{Eq:DihedralConstruction})$. Now,
$x*a=xa$ and $a*x=ax$ for every $x\in G$, $a\in F$. Since $F$ is central in
$G$, $(F,\,*)$ is also central in $(G,\,*)$. Finally, $x*x$ belongs to $F$ for
every $x\in G$, thus $(G,\,*)/(F,\,*)$ is an elementary abelian $2$-group. By
Theorems \ref{Th:Cyclic} and \ref{Th:Dihedral}, $(G,\,*)$ is a Moufang loop.
Then Theorem \ref{Th:Hsu} implies that $(G,\,*)$ is a code loop. Another
consequence of Theorems \ref{Th:Cyclic} and \ref{Th:Dihedral} is that
$N(G)=N(G,\,*)$. Hence the radical of the associator map $P_3$ in $G$ coincides
with the radical of the associator map $P^*_3$, where $P^*$ is the power map in
$(G,\,*)$.

To prove (iii), suppose that $\rad{P_3}=F$. Then $h\in N(G)\subseteq F\subseteq
Z(G)\subseteq N(G)$, where the only nontrivial inclusion $N(G)\subseteq F$
follows from the fact that $\rad{P_3}$ is trivial.
\end{proof}

Consider this general result about Moufang loops and code loops with trivial
radical.

\begin{lemma}\label{Lm:EquivalentCode}
Suppose that $L$ is a Moufang loop whose associator is equivalent to the
associator of a code loop $G$ with trivial radical. Then $L$ is a code loop
with trivial radical.
\end{lemma}
\begin{proof}
By the assumptions, $A(G)\le N(G)=Z(G)$, therefore $A(L)\le N(L)=Z(L)$, and
$L/N(L)$ is a group. Let $R$ be the associator map in $L$, and let $x$, $y$,
$z\in L$. Then $R(x,\,y,\,z)=0$ if and only if $R(x^{-1},\,y,\,z)=0$, by the
Moufang theorem. Since $|A(L)|\le 2$, we obtain
\begin{equation}\label{Eq:R}
    R(x,\,y,\,z)=R(x^{-1},\,y,\,z)
\end{equation}
for every $x$, $y$, $z\in L$. Because $R$ is equivalent to the associator map
of the code loop $G$, it is trilinear and $\rad{R}=N(L)$. Then $(\ref{Eq:R})$
implies $xN(L)=x^{-1}N(L)$ in $L/N(L)$, and $L/N(L)$ is an elementary abelian
$2$-group.
\end{proof}

\begin{lemma}\label{Lm:SuggestedByReferee}
Assume that $h\in F$, and that $(G,\,*)$ is constructed from a code loop
$G=(V,F,P)$ as in Lemma $\ref{Lm:Code}$. Let $P^*$ be the power map of
$(G,\,*)$. When $G/S\simeq C_2$ then
\begin{equation}\label{Eq:PowerMap}
    \begin{array}{c}
    P^*(xF)=\left\{\begin{array}{ll}
        P(xF),&x\in S,\\
        P(xF)+h,&x\in G\setminus S,
    \end{array}\right.
    \end{array}
\end{equation}
and $P^*-P$ is linear.

Else $G/S\simeq V_4$,
\begin{equation}\label{Eq:PowerMap2}
    \begin{array}{c}
    P^*(xF)=\left\{\begin{array}{ll}
        P(xF),&x\not\in \alpha,\\
        P(xF)+h,&x\in\alpha,
    \end{array}\right.
    \end{array}
\end{equation}
$($where $\alpha=\beta\gamma$ is as usual$)$, and $P^*-P$ is a quadratic form.
\end{lemma}
\begin{proof}
Since $x*y\in\{xy,\,xyh\}$, the addition in $G/F$ coincides with the addition
in $(G,\,*)/F$, and we can let $G/F=(G,\,*)/F=V$. By Lemma \ref{Lm:Code}(i),
$G/S\simeq C_2$ or $G/S\simeq V_4$. If $G/S\simeq C_2$, we have
$(\ref{Eq:PowerMap})$. Thus $P^*-P$ is linear.

If $G/S\simeq V_4$, we have $(\ref{Eq:PowerMap2})$. We claim that $R=P^*-P$ is
a quadratic form. First of all, $R_2(xF,\,yF)=R(xF)+R(yF)+R(xF+yF)$ does not
vanish if and only if $x$, $y$ belong to $\alpha\cup\beta\cup\gamma$ but not to
the same coset at the same time. Then $R_3(xF,\,yF,\,zF) = R_2(xF,\,zF) +
R_2(yF,\,zF) + R_2(xF+yF,\,zF)$ always vanishes, as one easily checks.
\end{proof}

We are ready to characterize all loops obtainable from code loops with trivial
radical via both of the constructions. We will also show how to connect all
code loops with the same associator maps.

\begin{proposition}\label{Pr:AllCode}
Let $G=(V,F,P)$ be a code loop with power map $P$. Let $H_0=G$, $H_1$, $\dots$,
$H_s$ be a sequence of loops, where $H_{i+1}$ is obtained from $H_i$ by the
cyclic or the dihedral construction, for $i=0$, $\dots$, $s-1$. If $\rad{P_3}$
is trivial, then $H_s$ is a code loop with power map $R$ satisfying
$\cdeg{(R-P)}\le 2$. Whether $\rad{P_3}$ is trivial or not, every code loop
$H_s$ with power map $R$ satisfying $\cdeg{(R-P)}\le 2$ can be obtained from
$H_0$ in this way.
\end{proposition}
\begin{proof}
Denote by $P^*$ the power map in $H_1$. For the rest of this paragraph, assume
that $P_3$ has trivial radical. By Lemma \ref{Lm:Code}, $H_1$ is a code loop
with trivial radical, and, by Lemma \ref{Lm:SuggestedByReferee},
$\cdeg{(P^*-P)}\le 2$. By induction, $H_s$ is a code loop and $\cdeg{(R-P)}\le
2$.

In fact, the two maps $P^*-P$ from $(\ref{Eq:PowerMap})$ and
$(\ref{Eq:PowerMap2})$ are available as long as $h\in F$, no matter what
$\rad{P_3}$ is.

In order to obtain all code loops with $\cdeg{(R-P)}\le 2$ from $H_0$, we must
show that the forms $P^*-P$ from $(\ref{Eq:PowerMap})$ and
$(\ref{Eq:PowerMap2})$ generate all forms with $\cdeg{}\le 2$, i.e., all
quadratic forms. Every quadratic form $Q$ determines an alternating bilinear
form $Q_2$, and when $Q_2=T_2$ for two quadratic forms $Q$, $T$, their
difference $Q-T$ is a linear form. We must therefore show how to obtain all
linear forms, and also all alternating bilinear forms as second derived forms
of maps stemming from $(\ref{Eq:PowerMap})$ and $(\ref{Eq:PowerMap2})$.

Note that the difference $P^*-P$ in $(\ref{Eq:PowerMap})$ determines a
hyperplane $S\cap V$ of $V$. Conversely, if $W\le V$ is a hyperplane, then
$W+F$ is a normal subloop of $V+F$. In this way, we obtain all linear forms.

In $(\ref{Eq:PowerMap2})$, $Q=P^*-P$ is a quadratic form such that
$\rad{Q_2}=S$ has codimension $2$ (since $|G/S|=4$). Moreover,
$Q_2(\gamma,\,\gamma)=Q_2(\beta,\,\beta)=0$, $Q_2(\beta,\,\gamma)\ne 0$, so
that $Q=U\oplus S$ for a hyperbolic plane $U=\langle x,\,y\rangle$,
$x\in\beta$, $y\in\gamma$. In this way, we can obtain all hyperbolic planes.
Every alternating bilinear form $f$ can be expressed as $U_1\oplus\cdots\oplus
U_k\oplus\rad{f}$, where every $U_i$ is a hyperbolic plane. Thus, by summing up
the differences $Q$ from repeated applications of the dihedral construction, we
can obtain any alternating bilinear form.
\end{proof}

Let us summarize the results about code loops obtained in this section:

\begin{theorem}\label{Th:Code}
If $G$ is a code loop with trivial radical and $\cc$ or $\dc$ is satisfied for
some $S\le G$, then $G/S$ is isomorphic to $C_2$ or $V_4$. The resulting loop
$(G,\,*)$ is a code loop with trivial radical, and the associators of $G$ and
$(G,\,*)$ are equivalent. Every Moufang loop whose associator is equivalent to
the associator of a code loop with trivial radical is itself a code loop with
trivial radical. Finally, any two code loops with equivalent associators can be
connected by the cyclic and dihedral constructions, possibly repeated.
\end{theorem}

\begin{remark}
It is not hard to check that trilinear alternating forms with trivial radical
exist in dimension $n$ if and only if $n=3$ or $n\ge 5$. (There are many
nonequivalent trilinear alternating forms with trivial radical when $n\ge 9$.)
Consequently, there are code loops with trivial radical (i.e., with two-element
nucleus) of order $2^n$ if and only if $n=4$ or $n\ge 6$.
\end{remark}

\section{Loops of Type $M(G,\,2)$}\label{Sc:Chein}

\noindent Chein \cite{Chein} discovered the following way of building up
nonassociative Moufang loops from nonabelian groups: Let $G$ be a finite group,
and denote by $\chn{G}$ the set of new elements $\{\chn{x};\;x\in G\}$. Then
$M(G,\,2)=(G\cup \chn{G},\,\circ)$ with multiplication $\circ$ defined by
\begin{equation}\label{Eq:Chein}
    x\circ y=xy,\dispskip x\circ\chn{y}=\chn{yx},\dispskip
    \chn{x}\circ y=\chn{xy^{-1}},\dispskip \chn{x}\circ\chn{y}=y^{-1}x
\end{equation}
is a Moufang loop that is associative if and only if $G$ is abelian. As the
restriction of the multiplication $\circ$ on $G$ coincides with the
multiplication in $G$, we will usually denote the multiplication in $M(G,\,2)$
by $\cdot$, too.

Many small Moufang loops are of this type; for instance $16/k$ for $k\le 2$,
and $32/k$ for $k\le 9$, where $n/k$ is the $k$th nonassociative Moufang loop
of order $n$. (See Section \ref{Sc:Small} for details. Table $1$ in \cite[p.\
A-3]{GMR} lists all loops $M(G,\,2)$ of order at most $63$.)

In this Section we are going to explore the effects of our constructions on
loops $M(G,\,2)$. The results are summarized in Corollary \ref{Cr:AllCyclic}
for the cyclic construction, and in Proposition \ref{Pr:AllDihedral} for the
dihedral construction.

The following Lemma gives some basic properties of loops $M(G,\,2)$:

\begin{lemma}\label{Lm:MG2Basic}
Let $G$ be a group and let $L=M(G,\,2)$ be the Moufang loop defined above.
Then:
\begin{enumerate}
\item[(i)] If $G$ is an abelian group then $N(L)=L$, else $N(L)=Z(G)$.

\item[(ii)] If $G$ is an elementary abelian $2$-group then $Z(L)=L$, else
$Z(L)=Z(G)\cap\{x\in G;\;x^2=1\}$.

\item[(iii)] If $S\le L$ then $S\le G$ or $|S\cap G|=|S\cap\chn{G}|$.

\item[(iv)] If $S\unlhd G$ then $S\unlhd L$.

\item[(v)] If $S\unlhd L$ then $S\unlhd G$, or both $G/(S\cap G)$ and $L/S$ are
elementary abelian $2$-groups.
\end{enumerate}
\end{lemma}
\begin{proof}
We know that $N(L)=L$ if and only if $G$ is abelian. Assume that $G$ is not
abelian. Then there are $x$, $y$, $z\in G$ such that $\chn{x}\cdot yz =
\chn{x(yz)^{-1}} \ne \chn{xy^{-1}z^{-1}} = \chn{x}y\cdot z$, and thus no
element of $\chn{G}$ belongs to $N(L)$. We have $x\cdot y\chn{z}=\chn{zyx}$,
while $xy\cdot \chn{z}=\chn{zxy}$. Also, $x(\chn{y}\cdot\chn{z})=xz^{-1}y$,
while $x\chn{y}\cdot\chn{z}=z^{-1}yx$. Hence $x\in G$ belongs to $N(L)$ if and
only if $x\in Z(G)$. This proves (i).

When $G$ is an elementary abelian $2$-group, we have $L\simeq G\times C_2$. As
$x\chn{y}=\chn{yx}$ and $\chn{y}x=\chn{yx^{-1}}$, an element $x\in G$ commutes
with all elements of $L$ if and only if $x\in Z(G)$ and $x^2=1$. This proves
(ii).

Part (iii) is an easy exercise (or see \cite[Proposition 4.5]{Thesis}).

Let $S\unlhd G$, and let $\varphi:G\to H$ be a group homomorphism with kernel
$S$. It is then easy to see that $\psi:M(G,\,2)\to M(H,\,2)$ defined by
$\psi(g)=\varphi(g)$, $\psi(\chn{g})=\chn{\varphi(g)}$, for $g\in G$, is a
homomorphism of Moufang loops with kernel $S$. Thus $S\unlhd M(G,\,2)$, and
(iv) is proved.

Finally, assume that $S\unlhd L$ and $S\not\le G$. Then there is $y\in G$ such
that $\chn{y}\in S$. For every $x\in G$, the element
$x\chn{y}x^{-1}\cdot\chn{y}$ belongs to $S$, since $S\unlhd L$. However,
$x\chn{y}x^{-1}\cdot\chn{y}=\chn{yxx}\cdot\chn{y} = y^{-1}yxx=xx$. That is why
$S\cap G$ contains all squares $x^2$, for $x\in G$, and the group $G/(S\cap G)$
must be an elementary abelian $2$-group. Also, $\chn{x}\cdot\chn{x}=1$ for
every $x\in G$. Hence $L/S$ is an elementary abelian $2$-group.
\end{proof}

We now investigate the two constructions for loops $M(G,\,2)$.

\begin{lemma}\label{Lm:MG2More}
Let $G$ be a group and let $L=M(G,\,2)$ be the Moufang loop defined above.
Then:
\begin{enumerate}
\item[(i)] If $(G,\,S,\,\alpha,\,h)$ satisfies $\cc$ then $L/S$ is dihedral,
$h\in N(L)$, and $hxh=x$ for every $x\in L\setminus G$.

\item[(ii)] If $L/S$ is cyclic then $L/S\simeq C_2$ and either $S=G$ or
$G/S\cap G\simeq C_2$.
\end{enumerate}
\end{lemma}
\begin{proof}
Assume that $S\unlhd G$ and $G/S=\langle\alpha\rangle$ is cyclic of order $m$.
Set $a=\alpha$, $b=\chn{S}=\chn{\alpha^0}$. Then $\langle a,\,b\rangle=L/S$
and, thanks to diassociativity, $L/S$ is a group. Moreover, $a^m=S$,
$b^2=\chn{S}\cdot\chn{S}=S$, and $aba = \alpha\chn{\alpha^0}\alpha =
\chn{\alpha}\alpha = \chn{\alpha^0}=b$. We know from Lemma \ref{Lm:MG2Basic}(i)
that $h\in S\cap Z(G)$ belongs to $N(L)$. Pick $\chn{g}\in\chn{G}$. Then
$h\chn{g}h=\chn{gh}h=\chn{ghh^{-1}}=\chn{g}$. This proves (i).

We proceed to prove (ii). Assume that $L/S=\langle\alpha\rangle$ is cyclic.
There must be some $x\in G$ such that $\chn{x}\in\alpha$, else $\alpha\subseteq
G$, which is impossible. As $\chn{x}\cdot\chn{x}=1$, we have $\alpha^2=S$, and
$L/S\simeq C_2$ follows. The rest is obvious.
\end{proof}

Consider this generalization of loops $M(G,\,2)$, also found in \cite[Theorem
2']{Chein}: Let $G$ be a group, $\theta$ an antiautomorphism of $G$, and $1\ne
h\in Z(G)$ such that $\theta$ is an involution, $\theta(h)=h$, and
$x\theta(x)\in Z(G)$ for every $x\in G$. Then the loop
$M(G,\,\theta,\,h)=(G\cup\chn{G},\,\circ)$ with multiplication $\circ$ defined
by
\begin{equation}\label{Eq:CheinH}
    x\circ y = xy,\dispskip x\circ\chn{y}=\chn{yx},\dispskip
    \chn{x}\circ y=\chn{x\theta(y)},\dispskip\chn{x}\circ\chn{y}
    =\theta(y)xh,
\end{equation}
is a Moufang loop that is associative if and only if $G$ is abelian.

Notice how the multiplication in $M(G,\, {}^{-1},\,h)$ differs from that of
$M(G,\,2)$ only at $\chn{G}\times\chn{G}$.

We claim that $M(G,\,{}^{-1},\,h)$ is never isomorphic to $M(H,2)$, for any
groups $G$, $H$: Every element of $\chn{H}$ in $M(H,\,2)$ is an involution.
Calculating in $M(G,\,{}^{-1},\,h)$, we get $\chn{x}*\chn{x}=h$ for every $x\in
G$. Thus every element of $\chn{G}$ in $M(G,\,{}^{-1},\,h)$ is of order $2|h|$,
where $|h|$ is the order of $h$. Then there are simply not enough elements of
order $2|h|$ in $M(H,\,2)$ for $M(H,\,2)$ to be isomorphic to
$M(G,\,{}^{-1},\,h)$.

Using Lemma \ref{Lm:MG2More} and the definitions $(\ref{Eq:Chein})$ and
$(\ref{Eq:CheinH})$, we get:

\begin{corollary}\label{Cr:AllCyclic}
Let $G$ be a group and let $L=M(G,\,2)$ be the Moufang loop defined above.
Assume that $(L,\,S,\,\alpha,\,h)$ satisfies $\cc$. Then $S=G$ or $G/(S\cap
G)\simeq C_2$. When $S=G$, the Moufang loop $(L,\,*)$ is isomorphic to
$M(G,\,{}^{-1},\,h)$. Every loop $M(G,\,{}^{-1},\,h)$ with $h^2=1$ can be
obtained in this way. When $G/(S\cap G)\simeq C_2$, then the multiplication in
$(L,\,*)$ is given by
\begin{equation}\label{Eq:MG2h}
    x*y=\left\{\begin{array}{ll}
        x\cdot y,&\text{if $x\in S$ or $y\in S$},\\
        (x\cdot y)h,&\text{otherwise},
    \end{array}\right.
\end{equation}
where $x$, $y\in L$, and where $\cdot$ is the multiplication in $L$.
\end{corollary}

With the classification \cite{GMR} available, one can often determine the
isomorphism type of $(L,\,*)$ from Corollary \ref{Cr:AllCyclic}. To illustrate
this point, assume that $(L=M(G,\,2),\,S,\,\alpha,\,h)$ satisfies $\cc$ and
that $S=G$. When $G=D_8$, the loop $L=M(D_8,\,2)$ contains $2$ elements of
order $4$. Hence $(L,\,*)$ must contain $2+8=10$ elements of order $4$, and it
turns out that the only such nonassociative Moufang loop of order $16$ is
$16/5$, according to \cite{GMR}. Similarly, $16/2=M(Q_8,\,2)$ always yields
$16/2$---the octonion loop of order $16$. If $L=24/1=M(D_{12},\,2)$, $(L,\,*)$
is isomorphic to $24/4$; if $L=32/9=M(Q_{16},\,2)$, $(L,\,*)$ is $32/38$, etc.

Now for the dihedral construction:

\begin{proposition}\label{Pr:AllDihedral}
Let $G$ be a group and let $L=M(G,\,2)$ be the Moufang loop defined above.
Assume that $(L,\,S,\,\beta,\,\gamma,\,h)$ satisfies $\dc$. Then $(L,\,*)$ is
isomorphic to $M(H,\,2)$ for some group $H$. Moreover, $S\unlhd G$, or
$L/S\simeq G/(S\cap G)\simeq V_4$. When $S\unlhd G$, then $(G,\,S,\,G\setminus
S,\,h)$ satisfies $\cc$, and the loop $(L,\,*)$ is equal to $M((G,\,*),\,2)$.
\end{proposition}
\begin{proof}
Assume that $(L,\,S,\,\beta,\,\gamma,\,h)$ satisfies $\dc$. Since the only
elementary abelian $2$-group that is also dihedral is $V_4$, Lemma
\ref{Lm:MG2Basic}(v) implies that $S\unlhd G$, or $L/S\simeq G/S\cap G\simeq
V_4$. When $S\unlhd G$, the group $G/S$ is obviously cyclic.

Suppose that $S\unlhd G$ and $\alpha=G\setminus S$. Then $(G,\,S,\,\alpha,\,h)$
satisfies $\cc$, and we can construct the group $(G,\,*)$. We are going to show
that the loop $(L,\,*)$ obtained from $L$ by the dihedral construction is equal
to $(L,\,\circ)=M((G,\,*),\,2)$, where we have denoted the operation by $\circ$
to avoid confusion.

Write $G=\bigcup_{i\in M}\alpha^i$. Without loss of generality, suppose that
$\chn{\alpha^i}=\alpha^i\gamma=\beta\alpha^{1-i}$ for every $i\in M$. Let
$x\in\alpha^i$ and $y\in\alpha^j$. We must show carefully that $x*y=x\circ y$,
$x*\chn{y}=x\circ\chn{y}$, $\chn{x}*y=\chn{x}\circ y$, and
$\chn{x}*\chn{y}=\chn{x}\circ\chn{y}$. Clearly, $x*y=x\circ y$. Also,
$x*\chn{y} = (x\cdot\chn{y})\cdot h^{-\sigma(i+j)} = \chn{yx}\cdot
h^{-\sigma(i+j)} = \chn{yxh^{\sigma(i+j)}} = \chn{y*x} = x\circ \chn{y}$.
Similarly, $\chn{x}*y = (\chn{x}\cdot y)\cdot h^{\sigma(1-i+j)} =
\chn{xy^{-1}}\cdot h^{\sigma(1-i+j)} = \chn{xy^{-1}h^{-\sigma(1-i+j)}} =
\chn{xy^{-1}h^{\sigma(i-j)}} = \chn{x*y^{-1}} = \chn{x}\circ y$, where we have
used the coset relation $\alpha^i\gamma = \beta\alpha^{1-i}$, and
$-\sigma(t)=\sigma(1-t)$. Finally, $\chn{x}*\chn{y} = (\chn{x}\cdot
\chn{y})\cdot h^{-\sigma(1-i+j)} = y^{-1}xh^{-\sigma(1-i+j)} =
y^{-1}xh^{\sigma(i-j)} = y^{-1}*x = \chn{x}\circ\chn{y}$.

It remains to show that $(L,\,*)=M(H,\,2)$ for some $H$ whenever $L/S$ is
dihedral. We take advantage of \cite[Theorem 0]{Chein}: If $Q$ is a
nonassociative Moufang loop such that every minimal generating set of $Q$
contains an involution, then $Q=M(H,\,2)$ for some group $H$.

Pick $x\in e\alpha^{1-i} = \alpha^if$. If $x\in G$ then $\alpha^2=S$, and $x*x
= x\cdot x=1$. If $x\not\in G$ then $x*x = x\cdot x\cdot h^{\sigma(1-i+i)}= 1$.
Because $\langle\alpha\rangle$ is a subloop of $(L,\,*)$, we have just shown
that every (minimal) generating set of $(L,\,*)$ contains an involution.
\end{proof}

We conclude this section with an example generalizing \cite{DOP}.

\begin{example} It is demonstrated in \cite{DOP} that $D_{2^n}$ can be obtained from
$Q_{2^n}$ via the cyclic construction, for $n>2$. Indeed, if $G=D_{2^n}=\langle
a,\,b\rangle$, then $\langle a\rangle = S\unlhd G$, $G/S\simeq C_2$,
$h=a^{2^{n-2}}\in Z(G)$, and $(G,\,S,\,a,\,h)$ satisfies $\cc$. The inverse of
$b$ in $(G,\,*)$ is $hb$, as $b*hb=bhbh=1$. Thus $a^{2^{n-1}}=1$,
$b*b=bbh=a^{2^{n-2}}$, $(b*a)*(a^{2^{n-2}}b) = ba*a^{2^{n-2}}b =
baa^{2^{n-2}}ba^{2^{n-2}}=bab=a^{-1}$, and $(G,\,*)\simeq Q_{2^n}$ follows.
Then, by Lemma \ref{Lm:MG2More}(ii), $L/S=M(D_{2^n},\,2)/S$ is dihedral of
order $4$, and $(L,\,S,\,\beta,\,\gamma,\,h)$ satisfies $\dc$, where we can
choose $\beta$, $\gamma$ so that $\alpha=\beta\gamma=G\setminus S$. Proposition
\ref{Pr:AllDihedral} then yields $(L,\,*) = M((G,\,*),\,2)\simeq
M(Q_{2^n},\,2)$.
\end{example}

\section{Small Moufang Loops}\label{Sc:Small}

\noindent Both the cyclic and dihedral constructions were studied for small
$2$-groups. In particular, using computers, the following question was answered
positively for groups of order $8$, $16$ and $32$ in \cite{Zhukavets}:
\emph{Given two groups $G$, $H$ of order $n$, is it possible to construct a
sequence of groups $G_0\simeq G$, $G_1$, $\dots$, $G_s\simeq H$ so that
$G_{i+1}$ is obtained from $G_i$ by means of the cyclic or the dihedral
construction?} The purpose of this section is to study an analogous question
for small Moufang loops, not necessarily of order $2^n$.

We will rely heavily on \cite{GMR}, where one finds multiplication tables of
all nonassociative Moufang loops of order less than $64$; one for each
isomorphism type. The book \cite{GMR} is based on Chein's classification
\cite{Chein}.

Following the notational conventions of \cite{GMR} closely, the $k$th Moufang
loop of order $n$ will be denoted by $n/k$. Whenever we refer to a
multiplication table of $n/k$, we always mean the one given in \cite{GMR}.

As we have mentioned in the Introduction, the only orders $n\le 32$ for which
there are at least two non-isomorphic nonassociative Moufang loops are $n=16$,
$24$, and $32$, with $5$, $5$, and $71$ loops, respectively.

For $n=24$ and $n=32$, all nonassociative Moufang loops of order $n$ can be
split into two subsets according to the size of their associator subloop (or
nucleus). Namely,
\begin{eqnarray*}
    A_{24}&=&\{24/1,\,24/3,\,24/4,\,24/5\},\\
    B_{24}&=&\{24/2\},\\
    A_{32}&=&\{32/1,\,\dots,\,32/6,\,32/10,\,\dots,32/26,\,32/29,\,32/30,\\
        &&32/35,\,32/36,\,32/39,\,\dots,\,32/71\},\\
    B_{32}&=&\{32/7,\,\dots,\,32/9,\,32/27,\,32/28,\,32/31,\,\dots,\,
        32/34,\,32/37,\,32/38\}.
\end{eqnarray*}
The size of the nucleus and the size of the associator subloop for loops in the
subsets $A_i$, $B_i$ are as follows:
\begin{displaymath}
    \begin{array}{c|c|c}
    \text{class}&\text{size of nucleus}&\text{size of associator subloop}\\
    \hline
    A_{24}&2&3\\
    B_{24}&1&4\\
    A_{32}&4&2\\
    B_{32}&2&4
    \end{array}
\end{displaymath}
All loops $16/k$, for $1\le k\le 5$, have associator subloop and nucleus of
cardinality $2$. Since the associator subloops do not change under our
constructions (cf.\ Theorems \ref{Th:Cyclic} and \ref{Th:Dihedral}), a loop
from set $A_i$ cannot be transformed to a loop from set $B_i$ via any of the
two constructions. The striking result is that the converse is also true:

\begin{theorem}\label{Th:Graph} For $n=16$, $24$, $32$, let $\mathcal G(n)$ be
a graph whose vertices are all isomorphism types of nonassociative Moufang
loops of order $n$, and where two vertices form an edge if a representative of
the second type can be obtained from a representative of the first type by one
of the two constructions. $($Lemmas $\ref{Lm:CycNormal}$ and
$\ref{Lm:DihNormal}$ guarantee that $\mathcal G(n)$ is not directed.$)$ Then:
\begin{enumerate}
\item[(i)] The graph $\mathcal G(16)$ is connected.

\item[(ii)] There are two connected components in $\mathcal G(24)$, namely
$A_{24}$ and $B_{24}$.

\item[(iii)] There are two connected components in $\mathcal G(32)$, namely
$A_{32}$ and $B_{32}$.
\end{enumerate}
In all cases, the connected components correspond to blocks of loops with
equivalent associator, and also to blocks of loops that have nucleus of the
same size.
\end{theorem}
\begin{proof}
The proof depends on machine computation that, together with detailed
information about exhaustive search for edges in $\mathcal G(n)$, will be
presented elsewhere. Our GAP libraries are available online \cite{GAP}.
\end{proof}

It is possible to select representatives of each connected component so that
they can be described in a uniform way. For instance, select representatives
$16/1=M(D_8,\,2)$, $24/1=M(D_{12},\,2)$, $24/2=M(A_4,\,2)$, $32/1=M(D_8\times
C_2,\,2)$, and $32/7=M(D_{16},\,2)$. See Section \ref{Sc:Chein} for the
definition of loops $M(G,\,2)$.

It is certainly of interest that, although the groups $D_{16}$ and $D_8\times
C_2$ are connected, the loops $M(D_{16},\,2)=32/7$ and $M(D_8\times
C_2,\,2)=32/1$ are not. This, in view of Proposition \ref{Pr:AllDihedral},
means that the groups $D_{16}$ and $D_8\times C_2$ cannot be connected via the
cyclic construction.

\begin{veryspecialcolumnspacing}
\begin{table}
\caption{Multiplication table of $32/1=M(D_8\times C_2,\,2)$.} \label{Tb:M321}
\begin{tiny}
\begin{displaymath}
\begin{array}{rrrrrrrr|rrrrrrrr|rrrrrrrr|rrrrrrrr}
  1 & 2 & 3 & 4 & 5 & 6 & 7 & 8 & 9 &10 &11 &12 &13 &14 &15 &16 &17 &18 &19 &20 &21 &22 &23 &24 &25 &26 &27 &28 &29 &30 &31 &32 \\
  2 & 3 & 4 & 1 & 6 & 7 & 8 & 5 &10 &11 &12 & 9 &14 &15 &16 &13 &18 &19 &20 &17 &24 &21 &22 &23 &26 &27 &28 &25 &32 &29 &30 &31 \\
  3 & 4 & 1 & 2 & 7 & 8 & 5 & 6 &11 &12 & 9 &10 &15 &16 &13 &14 &19 &20 &17 &18 &23 &24 &21 &22 &27 &28 &25 &26 &31 &32 &29 &30 \\
  4 & 1 & 2 & 3 & 8 & 5 & 6 & 7 &12 & 9 &10 &11 &16 &13 &14 &15 &20 &17 &18 &19 &22 &23 &24 &21 &28 &25 &26 &27 &30 &31 &32 &29 \\
  5 & 8 & 7 & 6 & 1 & 4 & 3 & 2 &13 &16 &15 &14 & 9 &12 &11 &10 &21 &22 &23 &24 &17 &18 &19 &20 &29 &30 &31 &32 &25 &26 &27 &28 \\
  6 & 5 & 8 & 7 & 2 & 1 & 4 & 3 &14 &13 &16 &15 &10 & 9 &12 &11 &22 &23 &24 &21 &20 &17 &18 &19 &30 &31 &32 &29 &28 &25 &26 &27 \\
  7 & 6 & 5 & 8 & 3 & 2 & 1 & 4 &15 &14 &13 &16 &11 &10 & 9 &12 &23 &24 &21 &22 &19 &20 &17 &18 &31 &32 &29 &30 &27 &28 &25 &26 \\
  8 & 7 & 6 & 5 & 4 & 3 & 2 & 1 &16 &15 &14 &13 &12 &11 &10 & 9 &24 &21 &22 &23 &18 &19 &20 &17 &32 &29 &30 &31 &26 &27 &28 &25 \\
  \hline
  9 &10 &11 &12 &13 &14 &15 &16 & 1 & 2 & 3 & 4 & 5 & 6 & 7 & 8 &25 &26 &27 &28 &29 &30 &31 &32 &17 &18 &19 &20 &21 &22 &23 &24 \\
 10 &11 &12 & 9 &14 &15 &16 &13 & 2 & 3 & 4 & 1 & 6 & 7 & 8 & 5 &26 &27 &28 &25 &32 &29 &30 &31 &18 &19 &20 &17 &24 &21 &22 &23 \\
 11 &12 & 9 &10 &15 &16 &13 &14 & 3 & 4 & 1 & 2 & 7 & 8 & 5 & 6 &27 &28 &25 &26 &31 &32 &29 &30 &19 &20 &17 &18 &23 &24 &21 &22 \\
 12 & 9 &10 &11 &16 &13 &14 &15 & 4 & 1 & 2 & 3 & 8 & 5 & 6 & 7 &28 &25 &26 &27 &30 &31 &32 &29 &20 &17 &18 &19 &22 &23 &24 &21 \\
 13 &16 &15 &14 & 9 &12 &11 &10 & 5 & 8 & 7 & 6 & 1 & 4 & 3 & 2 &29 &30 &31 &32 &25 &26 &27 &28 &21 &22 &23 &24 &17 &18 &19 &20 \\
 14 &13 &16 &15 &10 & 9 &12 &11 & 6 & 5 & 8 & 7 & 2 & 1 & 4 & 3 &30 &31 &32 &29 &28 &25 &26 &27 &22 &23 &24 &21 &20 &17 &18 &19 \\
 15 &14 &13 &16 &11 &10 & 9 &12 & 7 & 6 & 5 & 8 & 3 & 2 & 1 & 4 &31 &32 &29 &30 &27 &28 &25 &26 &23 &24 &21 &22 &19 &20 &17 &18 \\
 16 &15 &14 &13 &12 &11 &10 & 9 & 8 & 7 & 6 & 5 & 4 & 3 & 2 & 1 &32 &29 &30 &31 &26 &27 &28 &25 &24 &21 &22 &23 &18 &19 &20 &17 \\
 \hline
 17 &20 &19 &18 &21 &22 &23 &24 &25 &28 &27 &26 &29 &30 &31 &32 & 1 & 4 & 3 & 2 & 5 & 6 & 7 & 8 & 9 &12 &11 &10 &13 &14 &15 &16 \\
 18 &17 &20 &19 &22 &23 &24 &21 &26 &25 &28 &27 &30 &31 &32 &29 & 2 & 1 & 4 & 3 & 8 & 5 & 6 & 7 &10 & 9 &12 &11 &16 &13 &14 &15 \\
 19 &18 &17 &20 &23 &24 &21 &22 &27 &26 &25 &28 &31 &32 &29 &30 & 3 & 2 & 1 & 4 & 7 & 8 & 5 & 6 &11 &10 & 9 &12 &15 &16 &13 &14 \\
 20 &19 &18 &17 &24 &21 &22 &23 &28 &27 &26 &25 &32 &29 &30 &31 & 4 & 3 & 2 & 1 & 6 & 7 & 8 & 5 &12 &11 &10 & 9 &14 &15 &16 &13 \\
 21 &22 &23 &24 &17 &20 &19 &18 &29 &30 &31 &32 &25 &28 &27 &26 & 5 & 8 & 7 & 6 & 1 & 2 & 3 & 4 &13 &16 &15 &14 & 9 &10 &11 &12 \\
 22 &23 &24 &21 &18 &17 &20 &19 &30 &31 &32 &29 &26 &25 &28 &27 & 6 & 5 & 8 & 7 & 4 & 1 & 2 & 3 &14 &13 &16 &15 &12 & 9 &10 &11 \\
 23 &24 &21 &22 &19 &18 &17 &20 &31 &32 &29 &30 &27 &26 &25 &28 & 7 & 6 & 5 & 8 & 3 & 4 & 1 & 2 &15 &14 &13 &16 &11 &12 & 9 &10 \\
 24 &21 &22 &23 &20 &19 &18 &17 &32 &29 &30 &31 &28 &27 &26 &25 & 8 & 7 & 6 & 5 & 2 & 3 & 4 & 1 &16 &15 &14 &13 &10 &11 &12 & 9 \\
 \hline
 25 &28 &27 &26 &29 &30 &31 &32 &17 &20 &19 &18 &21 &22 &23 &24 & 9 &12 &11 &10 &13 &14 &15 &16 & 1 & 4 & 3 & 2 & 5 & 6 & 7 & 8 \\
 26 &25 &28 &27 &30 &31 &32 &29 &18 &17 &20 &19 &22 &23 &24 &21 &10 & 9 &12 &11 &16 &13 &14 &15 & 2 & 1 & 4 & 3 & 8 & 5 & 6 & 7 \\
 27 &26 &25 &28 &31 &32 &29 &30 &19 &18 &17 &20 &23 &24 &21 &22 &11 &10 & 9 &12 &15 &16 &13 &14 & 3 & 2 & 1 & 4 & 7 & 8 & 5 & 6 \\
 28 &27 &26 &25 &32 &29 &30 &31 &20 &19 &18 &17 &24 &21 &22 &23 &12 &11 &10 & 9 &14 &15 &16 &13 & 4 & 3 & 2 & 1 & 6 & 7 & 8 & 5 \\
 29 &30 &31 &32 &25 &28 &27 &26 &21 &22 &23 &24 &17 &20 &19 &18 &13 &16 &15 &14 & 9 &10 &11 &12 & 5 & 8 & 7 & 6 & 1 & 2 & 3 & 4 \\
 30 &31 &32 &29 &26 &25 &28 &27 &22 &23 &24 &21 &18 &17 &20 &19 &14 &13 &16 &15 &12 & 9 &10 &11 & 6 & 5 & 8 & 7 & 4 & 1 & 2 & 3 \\
 31 &32 &29 &30 &27 &26 &25 &28 &23 &24 &21 &22 &19 &18 &17 &20 &15 &14 &13 &16 &11 &12 & 9 &10 & 7 & 6 & 5 & 8 & 3 & 4 & 1 & 2 \\
 32 &29 &30 &31 &28 &27 &26 &25 &24 &21 &22 &23 &20 &19 &18 &17 &16 &15 &14 &13 &10 &11 &12 & 9 & 8 & 7 & 6 & 5 & 2 & 3 & 4 & 1
\end{array}
\end{displaymath}
\end{tiny}
\end{table}

\end{veryspecialcolumnspacing}

\begin{example}\label{Ex:Code}
Let us return to code loops. Their multiplication tables are easy to spot
thanks to this result of Chein and Goodaire \cite[Theorem 5]{CheinGoodaire}:
\emph{A loop $L$ is a code loop if and only if it is a Moufang loop with
$|L^2|\le 2$.} Here, $L^2$ denotes the set of all squares in $L$.

All loops $16/k$, $1\le k\le 5$, are code loops with trivial radical $($i.e.,
with nucleus of cardinality $2)$. In view of Proposition $\ref{Pr:AllCode}$ and
Theorem $\ref{Th:Graph}$, it suffices to establish this just for one loop
$16/k$; for example, the octonion loop of order $16$ is a code loop.

The loops $32/k$ are code loops for $k\in\{1$, $\dots$, $3$, $10$, $\dots$,
$22\}$, all with nontrivial radical. Markedly, it is possible to obtain a code
loop from a loop that is not code. Consider the loops $32/1=M(D_8\times
C_2,\,2)$ $($its multiplication table is given in Table $\ref{Tb:M321})$, and
the loop $32/4=M(16\Gamma_2c_1,\,2)$ (its multiplication table is given in
Table $\ref{Tb:ThreeEights})$. The group $16\Gamma_2c_1$ has presentation
$\langle a,\,b;\; a^4=b^4=(ab)^2=\commutator{a^2}{b}=1\rangle$. The loop $32/1$
is a code loop, while the loop $32/4$ is not, by the result of Chein and
Goodaire. They are connected, however, by Theorem $\ref{Th:Graph}$.
\end{example}

\section{Conjectures and Prospects}

\noindent Recall that given two Moufang loops (or groupoids) $(G,\,\circ)$,
$(G,\,*)$ defined on the same set $G$, their \emph{distance} $d(\circ,\,*)$ is
the cardinality of the set $\{(a,\,b)\in G\times G;\; a\circ b\ne a*b\}$.

\begin{veryspecialcolumnspacing}
\setlength{\fboxsep}{1pt}
\def\ff#1{\framebox{#1}}

\begin{table}
\caption{Multiplication tables of $32/4=M(16\Gamma_2c_1,\,2)$ and
$32/7=M(D_{16},\,2)$.} \label{Tb:ThreeEights}
\begin{tiny}
\begin{displaymath}
\begin{array}{rrrrrrrr|rrrrrrrr|rrrrrrrr|rrrrrrrr}
 1 &     2 &     3 &     4 &     5 &     6 &     7 &     8 &     9 &    10 &    11 &    12 &    13 &    14 &    15 &    16 &    17 &    18 &    19 &    20 &    21 &    22 &    23 &    24 &    25 &    26 &    27 &    28 &    29 &    30 &    31 &    32\\
 2 &     3 &     4 &\ff{ 1}&     6 &     7 &     8 &\ff{ 5}&    10 &    11 &    12 &\ff{ 9}&    14 &    15 &    16 &\ff{13}&    18 &    19 &    20 &\ff{17}&    22 &    23 &    24 &\ff{21}&\ff{32}&    29 &    30 &    31 &\ff{28}&    25 &    26 &    27\\
 3 &     4 &\ff{ 1}&\ff{ 2}&     7 &     8 &\ff{ 5}&\ff{ 6}&    11 &    12 &\ff{ 9}&\ff{10}&    15 &    16 &\ff{13}&\ff{14}&    19 &    20 &\ff{17}&\ff{18}&    23 &    24 &\ff{21}&\ff{22}&\ff{27}&\ff{28}&    25 &    26 &\ff{31}&\ff{32}&    29 &    30\\
 4 &\ff{ 1}&\ff{ 2}&\ff{ 3}&     8 &\ff{ 5}&\ff{ 6}&\ff{ 7}&    12 & \ff{9}&\ff{10}&\ff{11}&    16 &\ff{13}&\ff{14}&\ff{15}&    20 &\ff{17}&\ff{18}&\ff{19}&    24 &\ff{21}&\ff{22}&\ff{23}&\ff{30}&\ff{31}&\ff{32}&    29 &\ff{26}&\ff{27}&\ff{28}&    25\\
 5 &     6 &     7 &     8 &     1 &     2 &     3 &     4 &    13 &    14 &    15 &    16 &     9 &    10 &    11 &    12 &    21 &    22 &    23 &    24 &    17 &    18 &    19 &    20 &    29 &    30 &    31 &    32 &    25 &    26 &    27 &    28\\
 6 &     7 &     8 &\ff{ 5}&     2 &     3 &     4 &\ff{ 1}&    14 &    15 &    16 &\ff{13}&    10 &    11 &    12 &\ff{ 9}&    22 &    23 &    24 &\ff{21}&    18 &    19 &    20 &\ff{17}&\ff{28}&    25 &    26 &    27 &\ff{32}&    29 &    30 &    31\\
 7 &     8 &\ff{ 5}&\ff{ 6}&     3 &     4 &\ff{ 1}&\ff{ 2}&    15 &    16 &\ff{13}&\ff{14}&    11 &    12 &\ff{ 9}&\ff{10}&    23 &    24 &\ff{21}&\ff{22}&    19 &    20 &\ff{17}&\ff{18}&\ff{31}&\ff{32}&    29 &    30 &\ff{27}&\ff{28}&    25 &    26\\
 8 &\ff{ 5}&\ff{ 6}&\ff{ 7}&     4 &\ff{ 1}&\ff{ 2}&\ff{ 3}&    16 &\ff{13}&\ff{14}&\ff{15}&    12 &\ff{ 9}&\ff{10}&\ff{11}&    24 &\ff{21}&\ff{22}&\ff{23}&    20 &\ff{17}&\ff{18}&\ff{19}&\ff{26}&\ff{27}&\ff{28}&    25 &\ff{30}&\ff{31}&\ff{32}&    29\\
\hline
 9 &\ff{16}&\ff{11}&\ff{14}&    13 &\ff{12}&\ff{15}&\ff{10}&     1 &\ff{ 8}&\ff{ 3}&\ff{ 6}&     5 &\ff{ 4}&\ff{ 7}&\ff{ 2}&    25 &    26 &    27 &    28 &    29 &    30 &    31 &    32 &    17 &    18 &    19 &    20 &    21 &    22 &    23 &    24\\
10 &    13 &\ff{12}&\ff{15}&    14 &     9 &\ff{16}&\ff{11}&     2 &     5 &\ff{ 4}&\ff{ 7}&     6 &     1 &\ff{ 8}&\ff{ 3}&    26 &    27 &    28 &\ff{25}&    30 &    31 &    32 &\ff{29}&\ff{24}&    21 &    22 &    23 &\ff{20}&    17 &    18 &    19\\
11 &    14 &     9 &\ff{16}&    15 &    10 &    13 &\ff{12}&     3 &     6 &     1 &\ff{ 8}&     7 &     2 &     5 &\ff{ 4}&    27 &    28 &\ff{25}&\ff{26}&    31 &    32 &\ff{29}&\ff{30}&\ff{19}&\ff{20}&    17 &    18 &\ff{23}&\ff{24}&    21 &    22\\
12 &    15 &    10 &    13 &    16 &    11 &    14 &     9 &     4 &     7 &     2 &     5 &     8 &     3 &     6 &     1 &    28 &\ff{25}&\ff{26}&\ff{27}&    32 &\ff{29}&\ff{30}&\ff{31}&\ff{22}&\ff{23}&\ff{24}&    21 &\ff{18}&\ff{19}&\ff{20}&    17\\
13 &\ff{12}&\ff{15}&\ff{10}&     9 &\ff{16}&\ff{11}&\ff{14}&     5 &\ff{ 4}&\ff{ 7}&\ff{ 2}&     1 &\ff{ 8}&\ff{ 3}&\ff{ 6}&    29 &    30 &    31 &    32 &    25 &    26 &    27 &    28 &    21 &    22 &    23 &    24 &    17 &    18 &    19 &    20\\
14 &     9 &\ff{16}&\ff{11}&    10 &    13 &\ff{12}&\ff{15}&     6 &     1 &\ff{ 8}&\ff{ 3}&     2 &     5 &\ff{ 4}&\ff{ 7}&    30 &    31 &    32 &\ff{29}&    26 &    27 &    28 &\ff{25}&\ff{20}&    17 &    18 &    19 &\ff{24}&    21 &    22 &    23\\
15 &    10 &    13 &\ff{12}&    11 &    14 &     9 &\ff{16}&     7 &     2 &     5 &\ff{ 4}&     3 &     6 &     1 &\ff{ 8}&    31 &    32 &\ff{29}&\ff{30}&    27 &    28 &\ff{25}&\ff{26}&\ff{23}&\ff{24}&    21 &    22 &\ff{19}&\ff{20}&    17 &    18\\
16 &    11 &    14 &     9 &    12 &    15 &    10 &    13 &     8 &     3 &     6 &     1 &     4 &     7 &     2 &     5 &    32 &\ff{29}&\ff{30}&\ff{31}&    28 &\ff{25}&\ff{26}&\ff{27}&\ff{18}&\ff{19}&\ff{20}&    17 &\ff{22}&\ff{23}&\ff{24}&    21\\
\hline
17 &\ff{20}&\ff{19}&\ff{18}&    21 &\ff{24}&\ff{23}&\ff{22}&    25 &    30 &    27 &    32 &    29 &    26 &    31 &    28 &     1 &\ff{ 4}&\ff{ 3}&\ff{ 2}&     5 &\ff{ 8}&\ff{ 7}&\ff{ 6}&     9 &    14 &    11 &    16 &    13 &    10 &    15 &    12\\
18 &    17 &\ff{20}&\ff{19}&    22 &    21 &\ff{24}&\ff{23}&    26 &    31 &    28 &\ff{29}&    30 &    27 &    32 &\ff{25}&     2 &     1 &\ff{ 4}&\ff{ 3}&     6 &     5 &\ff{ 8}&\ff{ 7}&\ff{16}&     9 &    14 &    11 &\ff{12}&    13 &    10 &    15\\
19 &    18 &    17 &\ff{20}&    23 &    22 &    21 &\ff{24}&    27 &    32 &\ff{25}&\ff{30}&    31 &    28 &\ff{29}&\ff{26}&     3 &     2 &     1 &\ff{ 4}&     7 &     6 &     5 &\ff{ 8}&\ff{11}&\ff{16}&     9 &    14 &\ff{15}&\ff{12}&    13 &    10\\
20 &    19 &    18 &    17 &    24 &    23 &    22 &    21 &    28 &\ff{29}&\ff{26}&\ff{31}&    32 &\ff{25}&\ff{30}&\ff{27}&     4 &     3 &     2 &     1 &     8 &     7 &     6 &     5 &\ff{14}&\ff{11}&\ff{16}&     9 &\ff{10}&\ff{15}&\ff{12}&    13\\
21 &\ff{24}&\ff{23}&\ff{22}&    17 &\ff{20}&\ff{19}&\ff{18}&    29 &    26 &    31 &    28 &    25 &    30 &    27 &    32 &     5 &\ff{ 8}&\ff{ 7}&\ff{ 6}&     1 &\ff{ 4}&\ff{ 3}&\ff{ 2}&    13 &    10 &    15 &    12 &     9 &    14 &    11 &    16\\
22 &    21 &\ff{24}&\ff{23}&    18 &    17 &\ff{20}&\ff{19}&    30 &    27 &    32 &\ff{25}&    26 &    31 &    28 &\ff{29}&     6 &     5 &\ff{ 8}&\ff{ 7}&     2 &     1 &\ff{ 4}&\ff{ 3}&\ff{12}&    13 &    10 &    15 &\ff{16}&     9 &    14 &    11\\
23 &    22 &    21 &\ff{24}&    19 &    18 &    17 &\ff{20}&    31 &    28 &\ff{29}&\ff{26}&    27 &    32 &\ff{25}&\ff{30}&     7 &     6 &     5 &\ff{ 8}&     3 &     2 &     1 &\ff{ 4}&\ff{15}&\ff{12}&    13 &    10 &\ff{11}&\ff{16}&     9 &    14\\
24 &    23 &    22 &    21 &    20 &    19 &    18 &    17 &    32 &\ff{25}&\ff{30}&\ff{27}&    28 &\ff{29}&\ff{26}&\ff{31}&     8 &     7 &     6 &     5 &     4 &     3 &     2 &     1 &\ff{10}&\ff{15}&\ff{12}&    13 &\ff{14}&\ff{11}&\ff{16}&     9\\
\hline
25 &    30 &    27 &    32 &    29 &    26 &    31 &    28 &    17 &\ff{20}&\ff{19}&\ff{18}&    21 &\ff{24}&\ff{23}&\ff{22}&     9 &\ff{12}&\ff{11}&\ff{10}&    13 &\ff{16}&\ff{15}&\ff{14}&     1 &     6 &     3 &     8 &     5 &     2 &     7 &     4\\
26 &    31 &    28 &\ff{29}&    30 &    27 &    32 &\ff{25}&    18 &    17 &\ff{20}&\ff{19}&    22 &    21 &\ff{24}&\ff{23}&    10 &     9 &\ff{12}&\ff{11}&    14 &    13 &\ff{16}&\ff{15}&\ff{ 8}&     1 &     6 &     3 &\ff{ 4}&     5 &     2 &     7\\
27 &    32 &\ff{25}&\ff{30}&    31 &    28 &\ff{29}&\ff{26}&    19 &    18 &    17 &\ff{20}&    23 &    22 &    21 &\ff{24}&    11 &    10 &     9 &\ff{12}&    15 &    14 &    13 &\ff{16}&\ff{ 3}&\ff{ 8}&     1 &     6 &\ff{ 7}&\ff{ 4}&     5 &     2\\
28 &\ff{29}&\ff{26}&\ff{31}&    32 &\ff{25}&\ff{30}&\ff{27}&    20 &    19 &    18 &    17 &    24 &    23 &    22 &    21 &    12 &    11 &    10 &     9 &    16 &    15 &    14 &    13 &\ff{ 6}&\ff{ 3}&\ff{ 8}&     1 &\ff{ 2}&\ff{ 7}&\ff{ 4}&     5\\
29 &    26 &    31 &    28 &    25 &    30 &    27 &    32 &    21 &\ff{24}&\ff{23}&\ff{22}&    17 &\ff{20}&\ff{19}&\ff{18}&    13 &\ff{16}&\ff{15}&\ff{14}&     9 &\ff{12}&\ff{11}&\ff{10}&     5 &     2 &     7 &     4 &     1 &     6 &     3 &     8\\
30 &    27 &    32 &\ff{25}&    26 &    31 &    28 &\ff{29}&    22 &    21 &\ff{24}&\ff{23}&    18 &    17 &\ff{20}&\ff{19}&    14 &    13 &\ff{16}&\ff{15}&    10 &     9 &\ff{12}&\ff{11}&\ff{ 4}&     5 &     2 &     7 &\ff{ 8}&     1 &     6 &     3\\
31 &    28 &\ff{29}&\ff{26}&    27 &    32 &\ff{25}&\ff{30}&    23 &    22 &    21 &\ff{24}&    19 &    18 &    17 &\ff{20}&    15 &    14 &    13 &\ff{16}&    11 &    10 &     9 &\ff{12}&\ff{ 7}&\ff{ 4}&     5 &     2 &\ff{ 3}&\ff{ 8}&     1 &     6\\
32 &\ff{25}&\ff{30}&\ff{27}&    28 &\ff{29}&\ff{26}&\ff{31}&    24 &    23 &    22 &    21 &    20 &    19 &    18 &    17 &    16 &    15 &    14 &    13 &    12 &    11 &    10 &     9 &\ff{ 2}&\ff{ 7}&\ff{ 4}&     5 &\ff{ 6}&\ff{ 3}&\ff{ 8}&     1
\end{array}
\end{displaymath}
\end{tiny}
\end{table}

\end{veryspecialcolumnspacing}

Assume that $(G,\,*)$ is constructed from the Moufang loop $(G,\,\circ)$ via
one of the constructions. Then, as we hinted on in the title,
$d(\circ,\,*)=n^2/4$, where $n=|G|$. We conjecture that, similarly as for
groups, this is the smallest possible distance:

\begin{conjecture}\label{Cj:NonIso}
Every two Moufang $2$-loops of order $n$ in distance less than $n^2/4$ are
isomorphic.
\end{conjecture}

Since $A(G,\,*)=A(G,\,\circ)$ if $\cc$ or $\dc$ is satisfied, we wonder what is
the minimum distance of two Moufang loops with nonequivalent associator.

\begin{conjecture}\label{Cj:Nonequiv}
Two Moufang loops of order $n$ with nonequivalent associator are in distance at
least $3n^2/8$.
\end{conjecture}

This is illustrated in Table \ref{Tb:ThreeEights} for $n=32$, where one can
find multiplication tables of $32/4=M(16\Gamma_2c_1,\,2)$ and
$32/7=M(D_{16},\,2)$ the way they are listed in \cite{GMR}. To obtain the
multiplication table for 32/7, permute the $8\cdot 8=64$ framed triangular
regions by switching region $(2k, j)$ with region $(2k+1,j)$, for $k=0$,
$\dots$, $3$, $j=0$, $\dots$, $7$.

This does not mean that two loops with nonequivalent associator cannot be
closer. In fact, if a group multiplication table contains a subsquare
\begin{equation}\label{Eq:Subsquare}
\begin{array}{cc}a&b \\ b&a\end{array}
\end{equation}
and if the group is sufficiently large ($n\ge 6$), then the loop obtained by
switching $a$ and $b$ in $(\ref{Eq:Subsquare})$ cannot be associative.

We conclude the paper with a few suggestions for future research:
\begin{enumerate}
\item[1.] Decide whether two Moufang loops $M_0$, $M_s$ of order $n$ with
equivalent associator can be connected by a series of Moufang loops $M_0$,
$M_1$, $\dots$, $M_s$ so that the distance of $M_{i+1}$ from $M_i$ is $n^2/4$,
for $i=0$, $\dots$, $s-1$. (Note that additional constructions are needed
already for $n=64$.)

\item[2.] The main result of \cite{Dz} says that when the parameters of any of
the constructions are varied in a certain way, the isomorphism type of the
resulting group will not be affected. Can this be generalized to Moufang loops?
(See \cite{PVToward} for a step in this direction.)

\item[3.] Is there a general construction that preserves three quarters of the
multiplication table yet yields a Moufang loop with nonequivalent associator?

\item[4.] This paper attempts to launch a new approach to Moufang $2$-loops, by
obtaining them using group-theoretical constructions. One can envision a
similar programme for Bol loops modulo Moufang loops, for instance.

\item[5.] While this paper was under review, one of the authors has determined
by computer search that there are $4262$ nonassociative Moufang loops of order
$64$ that can be obtained from loops $M(G,\,2)$ by the two constructions, where
$G$ is a nonabelian group of order $32$. See \cite{PVToward} for more details.
Are there other nonassociative Moufang loops of order $64$?
\end{enumerate}

\section{Acknowledgement}

\noindent We would like to thank Edgar~G.~Goodaire for providing us with
electronic files containing multiplication tables of all nonassociative Moufang
loops of order at most $32$. We also thank the referee for many useful comments
that resulted in an improved exposition of the material.

\bibliographystyle{plain}

\end{document}


\title[Term Formulas for Two Constructions]
{Term Formulas for the Cyclic and Dihedral Constructions}
\author[P.~Vojt\v{e}chovsk\'y]{Petr Vojt\v{e}chovsk\'y}
\email{petr@iastate.edu}
\address{Department of Mathematics, Iowa State University, 400 Carver Hall,
Ames, IA, $50011$ USA} \urladdr{http://www.vojtechovsky.com}

\maketitle

\noindent This short note contains an alternative proof of the correctness of
the cyclic and dihedral constructions, as defined in \cite{AlesPetr}. Compare
this note with Theorems 6.3 and 7.3 of \cite{AlesPetr}. Our proof here is based
on identities only; it does not use extensions of Moufang loops by abelian
groups. The argument is rather technical and detailed, which is the reason why
we did not include it in \cite{AlesPetr}. Nevertheless it offers a different
view on the cyclic and dihedral constructions in \cite{AlesPetr}, and suggest
new generalizations.

This note is not self-contained, however, all notation, concepts and
conventions are as in \cite{AlesPetr}.

\section{More Modular Arithmetics}

\noindent Let $m$, $M$, $\oplus$, $\ominus$ and $\sigma$ be as in
\cite{AlesPetr}.

\begin{lemma}\label{Lm:Properties}
The following identities are satisfied by $\sigma$:
\begin{equation}\label{Eq:CyclicSigma}
    \begin{array}{l}
    \sigma(i+j)+\sigma(i\oplus j+i)+\sigma((i\oplus j)\oplus i+k)\\
        \phantom{mmmmi}
        = \sigma(i+k)+\sigma(j+i\oplus k)+\sigma(i+j\oplus(i\oplus k)),
    \end{array}
\end{equation}
\begin{equation}\label{Eq:010Sigma}
    \begin{array}{l}
    -\sigma(i+j)+\sigma(1-i\oplus j+i)+\sigma(1-j+k)\\
        \phantom{mmmmmmi}
        = \sigma(i+k)+\sigma(1-j+i\oplus k)-\sigma(i+j\ominus(i\oplus k)),
    \end{array}
\end{equation}
\begin{equation}\label{Eq:100Sigma}
    \begin{array}{l}
    \sigma(i,k)-\sigma(j+1-i\oplus k)
            -\sigma(i+1-(i\oplus k)\ominus j)\\
    \phantom{mmmi}
        =-\sigma(i+j)-\sigma(i\oplus j+1-i) + \sigma(\ominus{j}+k),
    \end{array}
\end{equation}
\begin{equation}\label{Eq:110Sigma}
    \begin{array}{l}
    \sigma(1-i+j)-\sigma(i\ominus j+i)+\sigma(1-(i\ominus j)\oplus i+k)\\
        \phantom{mmmmmmmmm}
        = \sigma(1-i+k)-\sigma(1-j+i\ominus k)
            +\sigma(1-i+j\ominus(i\ominus k)).
    \end{array}
\end{equation}
\end{lemma}
\begin{proof}
All identities are proved along similar lines. We always start with an identity
in $\oplus$, $\ominus$, and derive the corresponding identity in $\sigma$.
Sometimes we take advantage of
\begin{equation}\label{Eq:MinusSigma}
    \sigma(1-i)=-\sigma(i).
\end{equation}
The identities
\begin{eqnarray*}
    ((i\oplus j)\oplus i)\oplus k &=& i\oplus(j\oplus (i\oplus k)),\\
    ((i\oplus j)\ominus i)\ominus k&=& i\oplus(j\ominus(i\oplus k)),\\
    ((i\oplus k)\ominus j)\ominus i&=& k\oplus(i\ominus(i\oplus j)),\\
    ((i\ominus j)\oplus i)\ominus k &=& i\ominus(j\ominus(i\ominus k))
\end{eqnarray*}
eventually yield $(\ref{Eq:CyclicSigma})$, $(\ref{Eq:010Sigma})$,
$(\ref{Eq:100Sigma})$, and $(\ref{Eq:110Sigma})$, respectively. Formula
$(\ref{Eq:CyclicSigma})$ follows immediately. To see $(\ref{Eq:010Sigma})$,
simplify $(i\oplus j)\ominus i$ to $j$, and use $(\ref{Eq:MinusSigma})$. To see
$(\ref{Eq:100Sigma})$, simplify $i\ominus(i\oplus j)$ to $\ominus j$, and use
$(\ref{Eq:MinusSigma})$. To see $(\ref{Eq:110Sigma})$, use
$(\ref{Eq:MinusSigma})$.
\end{proof}

\section{The Cyclic Construction}\label{Sc:Cyclic}

\noindent Recall that the new multiplication $*$ in the cyclic case is
described by
\begin{equation}\label{Eq:CyclicDefinition}
    x*y = xyh^{\sigma(i+j)}.
\end{equation}

\begin{proposition}\label{Pr:CyclicConstruction}
Let $*$ be defined as above. Then $(G,\,*)$ is a Moufang loop. It is a group if
and only if $(G,\,.)$ is.
\end{proposition}
\begin{proof}
When $(G,\,\cdot)$ is associative, we can proceed as in \cite[Proposition
3.1]{AlesPetr}. We are going to show that the Moufang identity
$((x*y)*x)*z=x*(y*(x*z))$ is satisfied. Since $h\in Z(G)$, we have
\begin{eqnarray*}
    ((x*y)*x)*z &=& ((xy)x)z\cdot h^{\sigma(i+j)+\sigma(i\oplus j+i)+
        \sigma((i\oplus j)\oplus i+k)},\\
    x*(y*(xz)) &=& x(y(xz))\cdot h^{\sigma(i+k)+\sigma(j+i\oplus k) +
    \sigma(i+j\oplus(i\oplus k))}.
\end{eqnarray*}
We are done by $(\ref{Eq:CyclicSigma})$.
\end{proof}

\section{The Dihedral Construction}\label{Sc:Dihedral}

\noindent Recall that the multiplication $*$ in the case of dihedral
construction is given by
\begin{equation}\label{Eq:DihedralDefinition}
    x*y=(xy)h^{(-1)^r\sigma(i+j)}.
\end{equation}

Let us introduce the following notation. For $x_1$, $\dots$, $x_n\in G$, and
$\circ\in\{\cdot,\,*\}$, let
\begin{eqnarray*}
    L^\circ(x_1,\,\dots,\,x_n)& = &
        (\cdots((x_1\circ x_2)\circ x_3)\cdots)\circ x_n,\\
    R^\circ(x_1,\,\dots,\,x_n)& = &
        x_1\circ (x_2\circ (\cdots(x_{n-1}\circ x_n)\cdots)).
\end{eqnarray*}
Furthermore, write $X^\circ_{m,n}$ for $X^\circ(x_m,\,\dots,\,x_n)$,
$X^\circ_{m,n}e$ for $X^\circ(ex_m,\,x_{m+1},\,\dots,\,x_n)$ and
$X^\circ_{m,n}f$ for $X^\circ(x_m,\,\dots,\,x_{n-1},\,x_nf)$, where $X=L$, $R$.

\begin{lemma}\label{Lm:LeftRight}
Let $x_1$, $\dots$, $x_n\in G$. If $L_{1,n}=R_{1,n}$, $L^*_{1,n}=R^*_{1,n}$,
and $L_{1,n}e=R_{1,n}e$, then $L^*_{1,n}e=R^*_{1,n}e$. Similarly, if
$L_{1,n}=R_{1,n}$, $L^*_{1,n}=R^*_{1,n}$, and $L_{1,n}f=R_{1,n}f$, then
$L^*_{1,n}f=R^*_{1,n}f$.
\end{lemma}
\begin{proof}
We will prove only the second statement; the first being similar. Let $a$, $b$
be such that $L^*_{1,n}=L_{1,n}h^a$, $R^*_{1,n}=R_{1,n}h^b$. Since
$L^*_{1,n}=R^*_{1,n}$ and $L_{1,n}=R_{1,n}$, we may assume that $a=b$. Thus, if
we prove $L^*_{1,n}f = (L_{1,n}f)h^{-a}$ and $R^*_{1,n}f = (R_{1,n}f)h^{-b}$,
we will be done by $L^*_{1,n}f=(L_{1,n}f)h^{-a}=(R_{1,n}f)h^{-a} =
(R_{1,n}f)h^{-b} = R^*_{1,n}f$.

Let $r$, $s$, $\nu$ be such that $L_{1,n-1}\in\alpha^r\cup e\alpha^r$, and
$x_n\in(\alpha^s\cup e\alpha^s)\cap G_\nu$. Then there is $c$ such that
$L^*_{1,n}=L_{1,n-1}h^c*x_n = L_{1,n}h^{(-1)^\nu(c+\sigma(r+s))}$, and
$L^*_{1,n}f = L_{1,n-1}h^c*x_nf = L_{1,n}h^{(-1)^{\nu+1}(c+\sigma(r+s))}$.

Define $a_i$, $b_i$ by $R^*_{i,n}=R_{i,n}h^{a_i}$, $R^*_{i,n}f =
(R_{i,n}f)h^{b_i}$. Since $R^*_{n-1,n}f = (R_{n-1,n}f)h^{-a_n}$, we have
$a_n=-b_n$. Pick $i<n$, and assume that $a_j=-b_j$ for every $j$, $i<j\le n$.
Further assume that $x_i\in\alpha^u\cup e\alpha^u$, and
$R_{i+1,n}\in(\alpha^v\cup \alpha^vf)\cap G_\mu$. Then
$R^*_{i,n}=x_i*R_{i+1,n}h^{a_{i+1}} = R_{i,n}h^{(-1)^\mu\sigma(u+v) +
a_{i+1}}$, while $R^*_{i,n}f = x_i* (R_{i+1,n}f)h^{b_{i+1}} =
(R_{i,n}f)h^{(-1)^{\mu+1}\sigma(u+v)+b_{i+1}}$. Thus $a_i=-b_i$, and we are
through.
\end{proof}

\begin{proposition}\label{Pr:DihedralConstruction}
Let $*$ be defined as above. Then $(G,\,*)$ is a Moufang loop. It is a group if
and only if $(G,\,.)$ is.
\end{proposition}
\begin{proof}
When $(G,\,\cdot)$ is associative, we can proceed as in \cite[Proposition
4.3]{AlesPetr}. Whether $(G,\,\cdot)$ is associative or not, we are going to
verify $((x*y)*x)*z=x*(y*(x*z))$, for $x$, $y$, $z\in G$. Thanks to Lemma
\ref{Lm:LeftRight}, we can assume that $z\in\alpha^k$.

When $x$, $y\in G_0$, everything follows by Proposition
$\ref{Pr:CyclicConstruction}$. Let us have a look at the three remaining cases.

Assume that $x\in\alpha^i$, $y\in\alpha^jf$. Then $xy\in\alpha^{i\oplus
j}\gamma=\beta\alpha^{1-(i\oplus{j})}$, $(xy)x\in\alpha^{i\oplus
j}\gamma\alpha^i=\alpha^j\gamma=\beta\alpha^{1-j}$, $xz\in\alpha^{i\oplus k}$,
and $y(xz)\in\alpha^j\gamma\alpha^{i\oplus k}=\alpha^{j\ominus(i\oplus
k)}\gamma$. Therefore
\begin{eqnarray*}
    ((x*y)*x)*z &=&((xy)x)z\cdot h^{-\sigma(i+j)+\sigma(1-(i\oplus j)+i) +
        \sigma(1-j+k)},\\
    x*(y*(x*z)) &=&x(y(xz))\cdot h^{\sigma(i+k)+\sigma(1-j+i\oplus k)
        -\sigma(i+j\ominus(i\oplus k))}.
\end{eqnarray*}
We are done by $(\ref{Eq:010Sigma})$.

Assume that $x\in e\alpha^i$, $y\in\alpha^j$. Then $x\in\alpha^{1-i}\gamma$,
$xy\in\beta\alpha^{i\oplus j}$, $(xy)x\in\beta\alpha^{i\oplus
j}\beta\alpha^i=\alpha^{\ominus j}$, $xz\in\beta\alpha^{i\oplus
k}=\alpha^{1-(i\oplus k)}\gamma$, and $y(xz)\in\alpha^j\beta\alpha^{i\oplus
k}=\alpha^{1-((i\oplus k)\ominus j)}\gamma$. Therefore
\begin{eqnarray*}
    ((x*y)*x)*z &=&((xy)x)z\cdot h^{-\sigma(i+j)-\sigma(i\oplus j+1-i) +
        \sigma(\ominus j+k)},\\
    x*(y*(x*z)) &=&x(y(xz))\cdot h^{\sigma(i+k)-\sigma(j+1-(i\oplus k))
        -\sigma(i+1-(i\oplus k)\ominus j)}.
\end{eqnarray*}
We are done by $(\ref{Eq:100Sigma})$.

Finally assume that $x\in\alpha^if$, $y\in\alpha^j f$. Then
$x\in\beta\alpha^{1-i}$, $xy\in\alpha^i\gamma\alpha^j\gamma=\alpha^{i\ominus
j}$, $(xy)x\in\alpha^{(i\ominus j)\oplus i}\gamma=\beta\alpha^{1-((i\ominus
j)\oplus i)}$, $xz\in\alpha^i\gamma\alpha^k=\alpha^{i\ominus k}\gamma$,
$y(xz)\in\alpha^j\gamma\alpha^{i\ominus k}\gamma=\alpha^{j\ominus(i\ominus
k)}$. Therefore
\begin{eqnarray*}
    ((x*y)*x)*z &=&((xy)x)z\cdot h^{\sigma(1-i+j)-\sigma(i\ominus j+i) +
        \sigma(1-(i\ominus j)\oplus i+k)},\\
    x*(y*(x*z)) &=&x(y(xz))\cdot h^{\sigma(1-i+k)-\sigma(1-j+i\ominus k)
        +\sigma(1-i+j\ominus(i\ominus k))}.
\end{eqnarray*}
We are done by $(\ref{Eq:110Sigma})$.
\end{proof}

Note that the purpose of Lemma \ref{Lm:LeftRight} is solely to reduce the
number of cases to be considered in the proof of Proposition
\ref{Pr:DihedralConstruction}. Undoubtedly, one could dispose of Lemma
\ref{Lm:LeftRight} by introducing more identities as in Lemma
\ref{Lm:Properties}.

\bibliographystyle{plain}